\documentclass[12pt]{article}

\usepackage[latin1]{inputenc}
\usepackage{amsmath,amssymb}
\usepackage{amsmath,accents}
\usepackage{alltt,graphics,graphpap,color}
\usepackage[dvips]{graphicx}
\usepackage{amsfonts}
\usepackage{amscd}
\usepackage{mathrsfs}

\newtheorem{Theorem}{Theorem}[section]

\newtheorem{Proposition}[Theorem]{Proposition}
\newtheorem{Lemma}[Theorem]{Lemma}
\newtheorem{Corollary}[Theorem]{Corollary}

\newtheorem{Example}[Theorem]{Example}

\definecolor{ro}{rgb}{0,0,0}
\newcommand{\ro}{\color{ro}}

\newcommand{\rl}{\mathbb R }

\newcommand{\ep}{\varepsilon}

\newcommand{\mm}{{\cal M}}
\newcommand{\bmm}{{\overline {\cal M} }}

\newcommand{\nb}{\nabla}

\newcommand{\la}{\left\langle}
\newcommand{\ra}{\right\rangle}

\newcommand{\proof}{\emph{Proof. }}
\newcommand{\cvd}{\hfill$\square$ \bigskip}

\newcommand{\sss}{\mathbb S}
\newcommand{\qq}{\mathbb H}
\newcommand{\kk}{\mathbb K}
\newcommand{\pp}{\mathbb P}
\newcommand{\cc}{\mathbb C}
\newcommand{\cp}{\mathbb {CP}}

\newcommand{\gfs}{g_{FS}}
\newcommand{\tm}{T_{max}}

\newcommand{\hhhc}{\left| H\right|^6}
\newcommand{\aaa}{\left| A\right|^2}

\newcommand{\aaaq}{\left| A\right|^4}
\newcommand{\aao}{\accentset{\ \circ}{\left|A\right|}^2}
\newcommand{\aaoq}{\accentset{\ \circ}{\left|A\right|}^4}
\newcommand{\hhh}{\left| H\right|^2}
\newcommand{\hhhq}{\left| H\right|^4}
\newcommand{\auo}{\accentset{\circ\ }{\left|h_1\right|}^2}
\newcommand{\amoq}{\accentset{\circ\ }{\left|h_{-}\right|}^4}
\newcommand{\auoq}{\accentset{\circ\ }{\left|h_1\right|}^4}
\newcommand{\amo}{\accentset{\circ\ }{\left|h_{-}\right|}^2}

\newcommand{\nbh}{\left|\nabla H\right|^2}
\newcommand{\nba}{\left|\nabla A\right|^2}

\newcommand{\fs}{f_{\sigma}}

\newcommand{\rr}{\bar r}
\newcommand{\dt}{\displaystyle{\frac{\partial}{\partial t}}}

\begin{document}

\def\qed{\hbox{\hskip 6pt\vrule width6pt height7pt
depth1pt  \hskip1pt}\bigskip}


\title{Mean curvature flow of pinched submanifolds of \boldmath $\cp^n$}

\author{\sc G. Pipoli and C. Sinestrari}
\date{}

\maketitle

\begin{abstract}
We consider the evolution by mean curvature flow of a closed submanifold of the complex projective space. We show that, if the submanifold has small codimension and satisfies a suitable pinching condition on the second fundamental form, then the evolution has two possible behaviors: either the submanifold shrinks to a round point in finite time, or it converges smoothly to a totally geodesic limit in infinite time. The latter behavior is only possible if the dimension is even. These results generalize previous works by Huisken and Baker on the mean curvature flow of submanifolds of the sphere. 
\end{abstract}

\medskip

 \bigskip

\section{Introduction}\hspace{5 mm}  \setcounter{equation}{0}\noindent

Let $F_0:\mm \to \cp^n$ be a smooth immersion of a closed connected manifold in the complex projective space. We denote by $A$ the second fundamental form and  by $H$ the mean curvature vector associated with the immersion. 
The evolution of $\mm_0 = F_0(\mm)$ by mean curvature  flow is the one--parameter family of  immersions $F:\mm \times [0,\tm[ \,\to\cp^{n}$ satisfying
 \begin{equation}\label{1.1}
\left\{\begin{array}{l}
\dt F(p,t)= H,
\qquad p \in \mm, \, t \geq 0, \medskip \\
F(\cdot,0)=F_0.
\end{array}\right.
\end{equation}
We denote by $\mm_t=F(\mm,t)$ the evolution of $\mm_0$ at time $t$.
It is well known that this problem has a unique smooth solution up to some maximal time $\tm \leq \infty$. Moreover, if $\tm$ is finite the curvature of $\mm_t$ necessarily becomes unbounded as $t \to \tm$ and we say that the flow develops a singularity in finite time. The main theorem proved in this work is the following.

\begin{Theorem}\label{maincodim}
Let $\mm_0$ be a closed submanifold of $\cc\pp^n$ of real dimension $m$ and codimension $k=2n-m$.
Suppose either $n\geq 3$ and $k=1$, or $n \geq 7$ and $2\leq k<\frac{2n-3}{5}$ (equivalently, $2\leq k<\frac{m-3}{4}$). If at every point of $\mm_0$ the inequality
\begin{equation}\label{pinching_codim}
\aaa< \left\{\begin{array}{ll}
\displaystyle \frac{1}{m-1}\hhh+2 & \text{\rm if } k=1, \medskip\\
\displaystyle \frac{1}{m-1}\hhh+\frac{m-3-4k}{m} \quad & \text{\rm if } k\geq 2,
\end{array}\right.
\end{equation}
is satisfied, then the same holds on $\mm_t$ for all $0<t<\tm$. Moreover, one of the two following properties holds:
\begin{itemize}
\item [1)] $\tm<\infty$, and $\mm_t$ contracts to a point as $t \to \tm$,
\item [2)] $\tm=\infty$, and $\mm_t$ converges to a smooth totally geodesic submanifold as $t \to \tm$.
\end{itemize}
Case 2) can only occur if $m$ is even, and the limit submanifold is isometric to $\cc\pp^{\frac{m}{2}}$.
\end{Theorem}

An inequality of the form \eqref{pinching_codim} above is usually called a {\em pinching condition} on the second fundamental form. For instance, in the case $k=1$ it gives a bound on how much each principal curvature of the submanifold can differ from the others.

The above statement says in particular that in odd dimension a submanifold satisfying our assumptions necessarily shrinks to a point under mean curvature flow. We remark that this property is not proved directly: we show that the only alternative to a round point is the behavior in 2), but such a behavior is excluded for odd dimension because the only totally geodesic submanifolds of $\cc\pp^n$ with small codimension as in our hypotheses are isometric to a complex projective space.

When $\mm_t$ shrinks to a point in finite time as in case 1) above, one can show that, after an appropriate rescaling, it converges to an $m$-dimensional sphere, a behavior which is usually described as ``convergence to a round point'', see e.g. \cite[\S 9-10]{H1}, \cite[\S 6]{LXZ}. As a consequence, we obtain the following classification result.

\begin{Corollary}
Let $\mm_0$ satisfy the hypotheses of Theorem \ref{maincodim}. Then, if $m$ is odd, $\mm_0$ is diffeomorphic to an $\sss^{m}$, while if $k$ is even, $\mm_0$ is diffeomorphic either to an $\sss^{m}$ or to a $\cc\pp^{\frac m2}$. In every case $\mm_0$ is simply connected.
\end{Corollary}

The behavior of submanifolds evolving by mean curvature flow has been studied by several authors in the last decades, especially in the case of codimension one. The first fundamental result was obtained by Huisken \cite{H1}, who showed that any closed convex hypersurface in Euclidean space shrinks to a round point in finite time. He later proved \cite{H2} that the same holds for hypersurfaces in general Riemannian manifolds satisfying a stronger convexity condition which takes into account the geometry of the ambient space. A similar analysis has then been carried out by several authors for flows driven by speeds different from the mean curvature, and many convergence results to a round point are known for hypersurfaces satisfy suitable convexity requirements, see \cite{AMZ, AM} and the references therein. More recently, Andrews and Baker \cite{AB} have considered the mean curvature flow in the case of higher codimension, and proved the convergence to a round point for submanifolds of arbitrary codimension of the Euclidean space satisfying a suitable pinching condition. Similar results have then been obtained by Liu, Xu,Ye and Zhao for submanifolds of hyperbolic spaces \cite{LXYZ}  and of general Riemannian manifolds \cite{LXZ}.

By contrast, very few authors have considered cases where the mean curvature flow converges to a stationary limit. In the context of weak solutions, there is a quite general result by White \cite[Theorem 11.1]{W}, asserting that a mean convex solution either disappears in finite time or converges to a finite collection of stable minimal submanifolds. For classical solutions, results of this kind are known only in special cases. For the curve shortening flow, Grayson \cite{G} showed that an embedded curve in a Riemannian surface either shrinks to a round point or converges smoothly to a geodesic. When the dimension of the evolving submanifold is larger than one, other kinds of singularities can occur and an analogous statement can only be expected under suitable restrictions. Until now, a higher dimensional analogue of the results of \cite{G} has only been obtained for  submanifolds of the sphere, which have been studied by Huisken \cite{H3} for codimension one and by Baker \cite{Ba} for arbitrary codimension. The results in the two cases can be stated together as follows.

\begin{Theorem}\label{teorH3} {\rm \cite{H3,Ba}} Let $\mm_0$ be a closed $n$ dimensional submanifold of $\sss^{n+k}$, with $n \geq 2$, and suppose that we have on $\mm_0$
\begin{equation}\label{pinchHB}
\begin{array}{ll}
\aaa<\frac{1}{n-1}\hhh+2,\quad&\text{if $n\geq 4$, or $n=3$ and $k=1$}, \medskip \\
 \aaa<\frac{3}{4}\hhh+\frac{4}{3},\quad& \text{if $n=2$ and $k=1$}, \medskip \\
 \aaa<\frac{4}{3n}\hhh+\frac{2(n-1)}{3},\quad & \text{if $n=2,3$ and $k>1$}.
\end{array}
\end{equation}
Then one of the following holds:
\begin{itemize}
\item [1)] $\tm$ is finite and the $\mm_t$'s converge to a round point as $t\to\tm$,
\item [2)] $\tm$ is infinite and the $\mm_t$'s converge to a smooth totally geodesic hypersurface $\mm_{\infty}$, isometric to $\sss^n$.
\end{itemize}
\end{Theorem}

As underlined in the above statements, a key ingredient in all these results is the invariance under mean curvature flow of a pinching condition of the form $|A|^2 < a|H|^2+b$, for some $a>0$ and $b \in \rl$. The values of $a,b$ such that the invariance holds depend on the properties of the ambient manifold. If the ambient manifold is flat \cite{H1,AB}, or hyperbolic \cite{LXYZ}, or general \cite{LXZ}, the invariance can only be obtained for suitable values of $b \leq 0$, so that the condition rules out the possibility of a stationary limit. In the case of the sphere, it is possible instead to have invariance with some $b>0$, which allows the two possible behaviors of the above statements. In addition, a pinching condition with $b>0$ is substantially weaker: for example, in the case of codimension one it allows for some nonconvex hypersurfaces. Although the special structure of the sphere is used in an essential way in \cite{H3,Ba}, it is natural to expect that similar results should hold for more general ambient spaces of positive curvature.

The results of this paper confirm this expectation in the case of the complex projective space, showing that suitably pinched submanifolds evolving by mean curvature flow exhibit similar properties to the ones of the sphere. The complex projective space is a natural ambient space to consider beside the sphere, since it is a symmetric Einstein manifold with positive, but no longer constant, sectional curvature. The different structure of the Riemann curvature tensor complicates the study of the evolution of the curvature quantities with respect to \cite{H3,Ba}, and forces us to restrict our analysis to submanifolds with suitably small codimension. 

The paper is organized as follows. After recalling some notation and preliminary results, we prove in Section 3 the invariance of the pinching condition. In this part, in order to efficiently estimate the reaction terms in the evolution equations, it is crucial to choose normal and tangent frames strongly linked with the geometry of $\cc\pp^n$. In Section 4 we study the behavior of the norm of the traceless part of the second fundamental form, which is used to measure the improvement of pinching as the maximal time is approached. Since our estimates have additional lower order terms compared with \cite{H3,Ba}, to prove our main theorem we have to treat separately the cases of $\tm$ finite and $\tm$ infinite, which we do in Sections 5 and 6 respectively. The former case is more technically involved, and the convergence is obtained by integral estimates as in the previous papers, while for $\tm$ infinite the result follows from a more direct argument. Finally, in Section 7 we show that in the case of hypersurfaces our main result also holds for the mean curvature flow in quaternionic projective spaces.

\section{Preliminaries}  \setcounter{equation}{0}

The ambient manifold $\cp^n$ is a K\"ahler manifold of complex dimension $n$ with complex structure  $J$. It can be regarded as a real Riemannian manifold of dimension $2n$ endowed with the Fubini-Study metric $\gfs$. We denote the curvature tensor and the Levi-Civita connection of $(\cp^n,\gfs)$ with $\bar R$ and $\bar \nabla$ respectively. Then $\bar R$ has the explicit form, for all tangent vector fields $X,Y,Z,W,$
\begin{equation}\label{curvature_tensor}
\begin{array}{rcl}
\bar R(X,Y,Z,W) &=& \gfs(X,Z)\gfs(Y,W)-\gfs(X,W)\gfs(Y,Z)\\
& &+\gfs(X,JZ)\gfs(Y,JW)-\gfs(X,JW)\gfs(Y,JZ)\\
& &+2\gfs(X,JY)\gfs(Z,JW).
\end{array}
\end{equation}
In particular, the sectional curvature of a tangent plane spanned by two orthonormal vector fields $X$ and $Y$ is
\begin{equation}\label{sectional_curvature}
\bar K(X,Y)=1+3\gfs(X,JY)^2,
\end{equation}
therefore $1 \leq \bar K \leq 4$ and $\bar K =1$ (resp. $\bar K =4$) if and only if $JY$ is orthogonal (resp. tangent) to $X$. Moreover $(\cp^n, \gfs)$ is a symmetric space, so $\bar{\nb}\bar{R}=0$, and is an Einstein manifold with Einstein constant $2(n+1)$.

Let now $\mm$ be a closed submanifold of $\cp^n$, with induced metric $g$, curvature tensor $R$ and connection $\nabla$. The tangent and normal space to $\mm$ at a point $p$ are denoted by $T_p \mm$ and $N_p \mm$ respectively. Throughout the paper we denote by $m$ the dimension of $\mm$ and by $k=2n-m$ its codimension. Unless specified otherwise, Latin letters $i,j,l,...$ run from $1$ to $m$, Greek letters $\alpha, \beta, \gamma,...$ run from $m+1$ to $m+k$. 

Let $e_1,\dots,e_{m+k}$ be an orthonormal frame tangent to $\cp^n$ at a point of $\mm$, such that the first $m$ vectors are tangent to $\mm$ and the other ones are normal. With respect to this frame, the second fundamental form can be written 
$$
A=\sum_{\alpha}h^{\alpha}\otimes e_{\alpha},
$$
where the $h^{\alpha}=\left(h^{\alpha}_{ij}\right)$ are symmetric 2-tensors. The trace of the second fundamental form with respect to the metric $g$ is the mean curvature vector $H$:
$$
H=\sum_{\alpha}{\rm tr}\, h^{\alpha}e_{\alpha}=\sum_{\alpha}\sum_{ij}g^{ij}h^{\alpha}_{ij}e_{\alpha}.
$$
The traceless part of the second fundamental form is defined as $\displaystyle{\accentset{\circ}{A}=A-\frac{1}{m}H \otimes g}$, and its components are 
$\accentset{\circ}{h}^\alpha_{ij}=h^\alpha_{ij}-\frac{H}{m}^\alpha g_{ij}$, where $H^\alpha= \sum_{rs}g^{rs}h^{\alpha}_{rs}$.
In particular, the squared length satisfies $\displaystyle{\aao=\aaa-\frac{1}{m}\hhh}$.

If $\mm$ is a hypersurface, then the mean curvature vector is a multiple of the unit normal vector $\nu$ and satisfies
$$
H=-(\lambda_1+\dots+\lambda_m)\nu,
$$
where $\lambda_1\leq\cdots\leq\lambda_m$ are the principal curvatures. In addition, we have
$ \aaa=\lambda_1^2+\cdots+\lambda_m^2
$ and
\begin{equation}
\aao = \aaa-\frac{1}{m}\hhh=\frac{1}{m}\sum_{i<j}\left(\lambda_i-\lambda_j\right)^2,
\end{equation}
so that smallness of $\aao$ implies that the curvatures are close to each other.

The evolution equations of the main curvature quantities of a submanifold evolving by mean curvature flow in a general Riemannian space have been computed in \cite{AB} and \cite{Ba}. In our case, they take a simpler form because the ambient manifold is symmetric. We recall here the equations satisfied by $\hhh, \aaa$ and by the volume form $d\mu_t$ associated with the immersion at time $t$.

\begin{Lemma}\label{evoluzione_AH}
On a submanifold evolving by mean curvature flow in a symmetric ambient space we have
\begin{align*} 
&\nonumber 1)\quad \dt \hhh  =  \displaystyle{\Delta\hhh -2\left|\nb H\right|^2+2\sum_{i,j}\left(\sum_{\alpha}H^{\alpha}h_{ij}^{\alpha}\right)^2+2\sum_{l,\alpha,\beta}\bar R_{l\alpha l \beta}H^{\alpha}H^{\beta},}\\
& \nonumber
2) \quad \dt \aaa  =  \displaystyle{\Delta\aaa-2\left| \nb A\right|^2+2\sum_{\alpha, \beta}\left(\sum_{i,j}h_{ij}^{\alpha}h_{ij}^{\beta}\right)^2+2\sum_{i,j,\alpha,\beta}\left[\sum_p h_{ip}^{\alpha}h_{jp}^{\beta}-h_{ip}^{\beta}h_{jp}^{\alpha}\right]^2}\\
 &\nonumber
\phantom{3) \quad \dt \aaa  =}\displaystyle{+4\sum_{i,j,p,q}\bar R_{ipjq}\left(\sum_{\alpha}h_{pq}^{\alpha}h_{ij}^{\alpha}\right)-4\sum_{j,l,p}\bar R_{ljlp}\left(\sum_{i,\alpha}h_{pi}^{\alpha}h_{ij}^{\alpha}\right)}\\
 &\nonumber
\phantom{3) \quad \dt \aaa  =}\displaystyle{+2\sum_{l,\alpha,\beta}\bar R_{l\alpha l\beta}\left(\sum_{ij}h_{ij}^{\alpha}h_{ij}^{\beta}\right)-8\sum_{j,p,\alpha,\beta}\bar R_{jp\alpha\beta}\left(\sum_i h_{ip}^{\alpha}h_{ij}^{\beta}\right)},\\
&\nonumber 3)\quad 
\dt d\mu_t=-\hhh d\mu_t.
\end{align*}
\end{Lemma}

When the codimension is one these equations have a simpler form.

\begin{Lemma}\label{evoluzione_AH_hyp}
On a hypersurface evolving by mean curvature flow in a symmetric ambient space we have
\begin{align}
&\nonumber 1)\qquad \dt \hhh=\Delta \hhh-2\nbh+2\hhh\left(\aaa+\bar Ric(\nu,\nu)\right),\\
&\nonumber 2)\qquad\dt\aaa = \Delta\aaa-2\nba+2\aaa\left(\aaa+\bar Ric(\nu,\nu)\right)\\
&\nonumber\phantom{3)\qquad\dt\aaa =} -4 \sum_{i,j,p,l}\left(h_{ij}h_j^{\phantom j p}\bar R_{pli}^{\phantom{pli}l}-h^{ij}h^{lp}\bar R_{pilj}\right),
\end{align}
where $\bar Ric$ is the Ricci tensor of the ambient manifold.
\end{Lemma}


\section{Invariance of pinching}\setcounter{equation}{0}\setcounter{Theorem}{0}

In this section we prove that the pinching condition \eqref{pinching_codim} is invariant under the flow. To obtain the desired estimates, it is important to perform the computations using special tangent frames with suitable properties, which we now describe.

A first kind of frames, which was also considered in \cite{AB,LXZ}, can be defined at any point where $H\neq 0$ in the following way. We choose a privileged normal direction setting
\begin{equation}\label{base04}
e_{m+1}=\frac{H}{\left|H\right|}.
\end{equation}
Then we can choose $e_{m+2},\dots,e_{m+k}$ such that $\{e_{m+1},\dots,e_{m+k}\}$ is an orthonormal basis of $N_p\mm_t$  and choose any orthonormal basis $\{e_1,\dots,e_m\}$ of $T_p\mm_t$. Any tangent frame obtained in this way will be called of kind {\bf (B1)}.

With such a choice of tangent frame, the second fundamental form and its traceless part satisfy
$$
\left\{\begin{array}{ll}
\mbox{tr}\,h^{m+1}=\left|H\right|, \quad & \medskip \\
\mbox{tr} \, h^{\alpha}=0, &\alpha\geq m+2
\end{array}\right.
$$
and
$$
\left\{\begin{array}{ll}
\displaystyle{\accentset{\circ\phantom{^{nm1}}}{h^{m+1}}=h^{m+1}-\frac{\left|H\right|}{m}g, \quad} & \\
\accentset{\circ\phantom{^{\alpha}}}{h^{\alpha}}=h^{\alpha},&\alpha\geq m+2.
\end{array}\right.
$$
When using a basis of kind (B1), we adopt the following notation:
\begin{equation} \label{notazione_h}
\left|h_1\right|^2:=\left|h^{m+1}\right|^2,\qquad
\auo:=\accentset{\circ\phantom{^{m+1}}}{\left|h^{m+1}\right|}^2, \qquad
\left|h_-\right|^2=\amo:=\displaystyle{\sum_{\alpha=m+2}^{2n}\accentset{\circ\phantom{^{\alpha}}}{\left|h^{\alpha}\right|}^2.}
\end{equation}

A second kind of frames, more linked with the geometry of $\cc\pp^n$, is useful when we have to compute explicitly the components of the Riemann curvature tensor of the ambient manifold. The properties required in this case are described in the following lemma.

\begin{Lemma}\label{base_B2}
Let $\mm$ be a submanifold of $\cc\pp^n$ of dimension $m$ and codimension $k$. If $k \leq m$, then for every point $p\in\mm$ there exist $\{e_1,\dots,e_m\}$ an orthonormal basis of $T_p\mm$ and $\{e_{m+1},\dots,e_{m+k}\}$ an orthonormal basis of $N_p\mm$ such that:
\begin{enumerate}
\item for every $r\leq\frac k2$ we have \begin{equation}\label{base01}
\left\{\begin{array}{lcl}
Je_{m+2r-1} & = & \tau_r e_{2r-1}+\nu_re_{m+2r}, \medskip \\
Je_{m+2r} & = & \tau_r e_{2r}-\nu_re_{m+2r-1},
\end{array}\right.
\end{equation}
with $\tau_r,\nu_r\in [0,1]$ and $\tau_r^2+\nu_r^2=1$.
\item If $k$ is odd then $Je_{m+k}=e_k$.
\item The remaining vectors satisfy
\begin{equation}\label{base02}
Je_{k+1}=e_{k+2}, Je_{k+3}=e_{k+4},\dots,Je_{m-1}=e_m.
\end{equation}
\end{enumerate}
\end{Lemma}
\proof  For every point $p\in\mm$ the function
$$
\begin{array}{rccl}
\varphi: &N_p\mm\times N_p\mm & \rightarrow & \mathbb R\\
 & (X,Y) & \mapsto & \varphi(X,Y):=g(JX,Y)
\end{array}
$$
is a skew-symmetric bilinear form. It is a well-known fact that there is an orthonormal basis $\{e_{m+1}, \dots,e_{m+k}\}$ of $N_p\mm$ such that $\varphi$ is represented by the matrix
$$
M_{\varphi} = \left(\begin{array}{cccc}
\begin{array}{cc}
0 & \nu_1\\
-\nu_1 & 0
\end{array} 
& 0 & \cdots & 0\\
0 & \begin{array}{cc}
0 & \nu_2\\
-\nu_2 & 0
\end{array} & & 0\\
\vdots& & \ddots&\vdots\\
0 & 0 & \cdots & \begin{array}{cc}
0 & \nu_s\\
-\nu_s & 0
\end{array}
\end{array}\right)
\qquad \text{if } k=2s,$$ 

$$
M_{\varphi} = \left(\begin{array}{ccccc}
\begin{array}{cc}
0 & \nu_1\\
-\nu_1 & 0
\end{array} 
& 0 & \cdots & 0&0\\
0 & \begin{array}{cc}
0 & \nu_2\\
-\nu_2 & 0
\end{array} & & 0&0\\
\vdots& & \ddots&\vdots&\vdots\\
0 & 0 & \cdots & \begin{array}{cc}
0 & \nu_s\\
-\nu_s & 0
\end{array} & 0\\
0 & 0 & \cdots & 0 & 0
\end{array}\right)
\qquad \text{if } k=2s+1.
$$
Using the property that $|\varphi(X,Y)| \leq |X| |Y|$, we find that $|\nu_r| \leq 1$ for any $r$, and after possibly reversing signs we can have $\nu_r \in [0,1]$. 

Observe first that if $k$ is odd statement 2 follows easily.  When we consider the other vectors of the basis, the above construction implies that, for every $r\leq \frac k2$, the normal component of $Je_{m+2r-1}$ is given by $\nu_r e_{m+2r}$, while the normal component of $Je_{m+2r}$ is given by $-\nu_re_{m+2r-1}$. Now let us distinguish the cases $\nu_r<1$ and $\nu_r=1$. In the first case, we have
\begin{equation}\label{je}
\left\{\begin{array}{lcl}
Je_{m+2r-1} & = & \tau_r T_{2r-1}+\nu_re_{m+2r},\\
Je_{m+2r} & = & \hat \tau_r T_{2r}-\nu_re_{m+2r-1},
\end{array}\right.
\end{equation}
where the $T_i$ are unit vectors of $T_p\mm$ and $\tau_r,\hat\tau_r\in\mathbb R$. 
The above relations imply
$$
\tau_r^2+\nu_r^2=1= \hat \tau_r^2+ \nu_r^2,
$$
so, up to changing the sign of $T_{2r-1}$ and $T_{2r}$, we can obtain $\tau_r=\hat\tau_r \in (0,1]$.

If instead $\nu_r=1$, this means that $Je_{m+2r-1}$ coincides with $e_{m+2r}$. In this case, we choose $T_{2r-1}$ to be any unit tangent vector which is orthogonal to $T_1,\dots,T_{2r-2}$ and which is also orthogonal to $Je_{m+1},\dots,Je_{m+k}$. It is easy to see that such a vector exists because of the assumption $k \leq m$. We then define $T_{2r}=JT_{2r-1}$. By construction, $T_{2r}$ is a tangent unit vector orthogonal to $T_1,\dots,T_{2r-1}$.
Observe that equations \eqref{je} hold also in this case, with $\tau_r=\hat \tau_r=0$.

In general, since $\{e_{m+1},\dots,e_{m+k}\}$ is an orthonormal basis, from equations \eqref{je}, we have for any $i\neq j$
$$
g(T_i,T_j)=0.
$$
Then we define $e_i=T_i$ for $i=1,\dots,k$, and we complete the basis of $T_p\mm$ in an ortho\-normal way by choosing $e_{k+1},\dots, e_m$ in such a way that requirement 3 is satisfied. \cvd

Any basis satisfying the properties of the previous lemma will be called of kind {\bf (B2)}. Since $J^2=-id$, from \eqref{base01} it follows easily that such a basis also satisfies
\begin{equation}\label{base03}
\left\{\begin{array}{lcl}
Je_{2r-1} & = & -\nu_r e_{2r}-\tau_re_{m+2r-1}, \medskip \\
Je_{2r} & = & \nu_r e_{2r-1}-\tau_re_{m+2r}.
\end{array}\right.
\end{equation}
If $k$ is odd, it is convenient to define $\tau_r=1$, $\nu_r=0$ for $r=\frac{k+1}{2}$. In this way, the first equations in \eqref{base01} and in \eqref{base03} hold also for this value of $r$.

In general, the requirements for (B1) and (B2) are incompatible and the two kinds of bases are different. Thus when we use frames of type B2, we have
$
H=\sum_{\alpha}H^{\alpha}e_{\alpha}, 
$ with $H^\alpha$ not necessarily zero for $\alpha>m+1$.

Observe that when $k=1$ these constructions are trivial: there is an unique (up to sign) normal unit vector $e_{2n}$, $H$ is a multiple of such vector and $e_1=Je_{2n}$ is a tangent vector. Then for a hypersurface we can choose a basis that is at the same time of type (B1) and (B2). \bigskip

When $k\geq 2$, we introduce the following notation taken from \cite{AB}
$$
R_1:=\sum_{\alpha, \beta}\left(\sum_{i,j}h_{ij}^{\alpha}h_{ij}^{\beta}\right)^2+\sum_{i,j,\alpha,\beta}\left[\sum_p h_{ip}^{\alpha}h_{jp}^{\beta}-h_{ip}^{\beta}h_{jp}^{\alpha}\right]^2,
$$
$$
R_2:=\sum_{i,j}\left(\sum_{\alpha}H^{\alpha}h_{ij}^{\alpha}\right)^2.
$$
\noindent If we use a frame of kind (B1), it is easily checked that
\begin{equation}\label{e.R2}
R_2=\left\{ \begin{array}{ll}
\displaystyle{\auo |H|^2 + \frac{1}{m}|H|^4 \qquad} & \mbox{if $H \neq 0$}  \medskip \\
0 & \mbox{if $H = 0$.}
\end{array}
\right. \end{equation}

The following result, proved in \cite[\S 3]{AB} and in \cite[\S 5.2]{Ba}, is useful in the estimation of the reaction terms occurring  in the evolution equations of Lemma \ref{evoluzione_AH}. It only uses the algebraic properties of $R_1$ and $R_2$ and is independent on the flow.

\begin{Lemma}\label{l.est}
At a point where $H \neq 0$ we have, for any $a \in \rl$
\begin{eqnarray*}
2R_1-2aR_2 & \leq & 2\auoq-2\left(a-\frac 2m \right) \auo |H|^2 - \frac 2m \left(a-\frac 1m \right)|H|^4 \\
& & + 8 \auo \amo + 3 \amoq. 
\end{eqnarray*}
In addition, if $a > 1/m$ and if $b \in \rl$ is such that $|A|^2=a|H|^2+b$, we have
\begin{eqnarray*}
2R_1-2aR_2 & \leq &  \left(6 - \frac{2}{ma-1}  \right) \aao\amo  -3\amoq \\
& & + \frac{2mab}{ma-1} \auo + \frac{4b}{ma-1} \amo - \frac{2b^2}{ma-1}.
\end{eqnarray*}

\end{Lemma}

We now derive a sharp estimate on the gradient terms appearing in the evolution equations for $\aaa$ and $\hhh$ which will be used many times in the rest of the paper. Observe that the results are independent of the flow. Our starting point is the following inequality, first proved in Lemma 2.2 of \cite{H2} in the case of hypersurfaces, and then extended to general codimension in Lemma 3.2 of \cite{LXZ}.

\begin{Lemma}\label{grad01}
Let $\bmm$ an Riemannian manifold and $\mm$ a submanifold of $\bmm$ of dimension $m$ and arbitrary codimension. Then
\begin{equation}\label{gradgen}
\nba \geq \left(\frac{3}{m+2}-\eta\right)\nbh-\frac{2}{m+2}\left(\frac{2}{m+2}\eta^{-1}-\frac{m}{m-1}\right)\left|\omega\right|^2,
\end{equation}
holds for any $\eta>0$. Here $\omega=\sum_{ij\alpha}\bar R_{\alpha jij}e_i\otimes \omega_{\alpha}$, where $\omega_\alpha$ is the dual frame to $e_\alpha$.
\end{Lemma}

Note that if the ambient space is Einstein, like $\cc\pp^n$, and if $\mm$ is a hypersurface, then $\omega=0$. So we can let $\eta \to 0$ in inequality \eqref{gradgen} and find

\begin{equation}\label{grad02} 
\nba\geq\frac{3}{m+2}\nbh.
\end{equation}

For submanifolds of higher codimension, $\omega$ is in general nonzero. However, using the special properties of $\cc\pp^n$, we can prove the following estimate.

\begin{Lemma}\label{grad03}
Let $\mm$ be a submanifold of $\cc\pp^n$ of dimension $m$ and codimension $k \leq m$. Then we have, at any point of $\mm$,
$$
\nba\geq \frac29(m+1) \left|\omega\right|^2.
$$
\end{Lemma}
\proof
We first compute explicitly $\left|\omega\right|^2$ using a basis of type (B2). The relations \eqref{base03} and the expression of $\bar R$ give
\begin{eqnarray*}
\bar R_{\alpha jij} & = & 3\gfs(e_{\alpha},Je_j)\gfs(e_i,Je_j) \\
&=& \left\{\begin{array}{ll}
3\tau_r\nu_r & \mbox{ if }\alpha=m+2r-1, \ i=2r,  \ j=2r-1, \medskip \\
-3\tau_r\nu_r \qquad & \mbox{ if } \alpha=m+2r, \ i=2r-1,\  j=2r, \medskip \\
0 & \mbox{ otherwise.}
\end{array}\right.
\end{eqnarray*}
We recall that if $k$ is odd then $\nu_{r}=0$ for $r=\frac{k+1}{2}$. Thus we have, for a general $k$,
\begin{equation}\label{omega}
\left|\omega\right|^2=18\sum_{r \leq \frac{k}{2}}\tau_r^2\nu_r^2.
\end{equation}

Next we recall a lower bound on $|\nabla A|$ for general submanifolds $\mm$ of $\cc\pp^n$ which was proved in \cite{Ko}. Following the notation of that paper, for any vector field $X$ tangent to $\mm$, we write $JX=PX+FX$, where $PX$ and $FX$ are the tangent and normal component of $JX$ respectively. Similarly, for a normal vector field $V$ we write  $JV=tV+fV$ where $tV$ is tangent to $\mm$ and $fV$ is normal. Then 
Lemma 3.6 of \cite{Ko} asserts that, at any point of $\mm$, we have
\begin{equation}\label{Kon}
\nba\geq 2\left(\left|P\right|^2\left|t\right|^2+\left|FP\right|^2\right). \medskip
\end{equation}
In a given orthonormal basis, the above norms are
$$\left|P\right|^2=\sum_{i=1}^{m}\left|Pe_i\right|^2,\quad \left|t\right|^2=\sum_{\alpha=m+1}^{2n}\left|te_{\alpha}\right|^2\ \text{and}\  \left|FP\right|^2=\sum_{i=1}^{m}\left|FPe_i\right|^2.$$
We choose again a basis of type (B2) and estimate the above expressions in the cases $k$ even and $k$ odd separately, using the relations \eqref{base01}, \eqref{base03}. If $k=2d$ we have
$$
\left|P\right|^2=(m-k)+\big( 2 \sum_{r \leq d} \nu_r^2  \big) = m- \big( 2 \sum_{r\leq d} \tau_r^2 \big),
$$
$$
\left|t\right|^2=2\sum_{r \leq d} \tau_r^2.
$$
Therefore, using the property $\nu_r^2+\tau_r^2=1$ and the assumption $m \geq k$, we find
\begin{eqnarray}
|P|^2|t|^2 & = & 2m \sum_{r \leq d} \tau_r^2 - 4 \sum_{r,s \leq d}\tau_r^2\tau_s^2 \nonumber \\
& \geq & 2m \sum_{r \leq d} \tau_r^2 - 2 \sum_{r,s \leq d}(\tau_r^4 + \tau_s^4) \nonumber \\
& = & 2m \sum_{r \leq d} \tau_r^2 - 2 k \sum_{r \leq d}\tau_r^4 \nonumber \\
& \geq & 2m \sum_{r \leq d}( \tau_r^2 - \tau_r^4 ) = \frac {m}{9} |\omega|^2. \label{evenPT}
\end{eqnarray}
If instead $k=2d+1$ we find
$$
\left|P\right|^2=(m-k)+\big( 2 \sum_{r \leq d} \nu_r^2  \big) = m-1- \big( 2 \sum_{r\leq d} \tau_r^2 \big),
$$
$$
\left|t\right|^2=1+2\sum_{r \leq d} \tau_r^2.
$$
Therefore,
\begin{eqnarray*}
|P|^2 |t|^2 & \geq & m-1+2(m-2)\sum_{r\leq d}\tau_r^2-2(k-1)\sum_{r\leq d}\tau_r^4\\
& \geq & m-1+2(m-2)\sum_{r\leq d}\tau_r^2\nu_r^2.
\end{eqnarray*}
Since for every $r$ we have $\nu_r^2+\tau_r^2=1$, we deduce that $\nu_r^2\tau_r^2\leq\frac 14$. Therefore, using that
$m-1 \geq k-1=2d$, we find
\begin{eqnarray}
|P|^2|t|^2 & \geq & 2d+2(m-2)\sum_{r\leq d}\tau_r^2\nu_r^2 \nonumber \\
& \geq & 2(m+2)\sum_{r\leq d}\tau_r^2\nu_r^2 = \frac{m+2}{9}|\omega|^2. \label{oddPT}
\end{eqnarray}

Finally, we have for any $k$
\begin{equation}
\left|FP\right|^2=2\sum_{r \leq \frac{k}{2}} \tau_r^2\nu_r^2 = \frac{|\omega|^2}{9}.\label{FP}
\end{equation}
Putting together inequalities \eqref{Kon}, \eqref{evenPT}, \eqref{oddPT} and \eqref{FP} the conclusion follows.
\cvd

The previous result allows us to obtain an estimate similar to \eqref{grad02} for general codimension.

\begin{Lemma}\label{grad05}
For any submanifold $\mm$ of $\cc\pp^n$ with dimension satisfying the assumptions of Theorem \ref{maincodim},
we have
$$
\nba \geq \frac{16}{9(m+2)} \nbh.
$$
\end{Lemma}
\proof If the codimension is $1$, then the result follows directly from \eqref{grad02}. In the case of higher codimension, the trick is to combine the estimates from Lemma \ref{grad01} and Lemma \ref{grad03} as follows: 

\begin{eqnarray*}
3 \nba  & = & 2\nba +  \nba \\
& \geq & 
2 \left(\frac{3}{m+2}-\eta \right)\nbh +\left[ \frac29(m+1)-\frac{4}{m+2}\left(\frac{2}{m+2}\eta^{-1}-\frac{m}{m-1}\right)\right]\left|\omega\right|^2.
\end{eqnarray*}
Now we choose $\eta=1/3(m+2)$ to obtain
\begin{eqnarray*}
 3\nba 
& \geq & 
\frac {16}3\frac{1}{m+2} \nbh + \left[ \frac29(m+1)-\frac{24}{m+2} \right]\left|\omega\right|^2,
\end{eqnarray*}
and the term inside square brackets is positive for $m$ as in our hypotheses. \cvd

We are now ready to prove the invariance of the pinching condition of Theorem \ref{maincodim}. We treat separately the case of hypersurfaces, where the analysis is simpler, and the case of higher codimension, where the two kinds of bases introduced before are essential. However, the strategy of proof is the same in the two cases: we consider
the function
$$Q=\aaa-a\hhh-b$$
for suitable $a,b$, 
and we analyze its evolution equation showing that, if $Q(x,t)=0$ at some point $(x,t)\in\mm\times\left[0,\tm\right[$, then $\left(\dt-\Delta\right)Q \leq 0$ at this point. By the maximum principle, the result will follow.

\begin{Proposition}\label{pinching_preserved_hyp}
Let $\mm_0$ be a closed hypersurface of $\cc\pp^n$, with $n\geq 3$. Then the pinching condition
\begin{equation}\label{pinching_stretto_hyp}
\aaa \leq \frac{1}{m-1+\ep} \hhh+2(1-\ep)
\end{equation}
is preserved by the mean curvature flow for any $\ep \in [0,1)$.
\end{Proposition}
\proof
Let us set $Q=\aaa-a\hhh-b$ with $a=(m-1+\ep)^{-1}$ and $b=2(1-\ep)$. Lemma \ref{evoluzione_AH_hyp} gives
\begin{equation}\label{evoluzione_Q_hyp}
\begin{array}{rcl}
\dt Q & =&  \Delta Q-2\left(\nba-a\nbh\right)+2\left(\aaa-a\hhh\right)\left(\aaa+\rr\right)\\
 & & -4\left(h_{ij}h_j^{\phantom j p}\bar R_{pli}^{\phantom{pli}l}-h^{ij}h^{lp}\bar R_{pilj}\right)\\
 & = & \Delta Q-2\left(\nba-a\nbh\right)+2Q\left(\aaa+\rr\right)+2b\left(\aaa+\rr\right)\\
 & & -4\left(h_{ij}h_j^{\phantom j p}\bar R_{pli}^{\phantom{pli}l}-h^{ij}h^{lp}\bar R_{pilj}\right),
\end{array}
\end{equation}
where we have set
\begin{equation}\label{barerre}
\rr=\bar Ric(\nu,\nu)=2(n+1).
\end{equation}
By Lemma \ref{grad05} the gradient terms in equation \eqref{evoluzione_Q_hyp} are non-positive and it suffices to consider the contribution of the reaction terms. Fix an orthonormal basis tangent to $\mm_t$ that diagonalizes the second fundamental form and call $\lambda_1\leq\lambda_2\leq \dots \leq\lambda_m$ its eigenvalues. Recalling that any sectional curvature $\bar K_{ij}$ satisfies $\bar K_{ij}\geq 1$, we find 
\begin{eqnarray}
\nonumber-4\left(h_{ij}h_j^{\phantom j p}\bar R_{pli}^{\phantom{pli}l}-h^{ij}h^{lp}\bar R_{pilj}\right) & = & -4\left(\lambda_j^2\delta_{ij}\delta_{jp}\bar R_{plil}-\lambda_j\lambda_l\delta_{ij}\delta_{lp}\bar R_{pilj}\right)\\
 \nonumber& = & \displaystyle{-4\sum_{j,l}\left(\lambda_j^2-\lambda_j\lambda_l\right)\bar R_{jljl}}\\
 \nonumber& = & \displaystyle{-2\sum_{j,l}\left(\lambda_j-\lambda_l\right)^2\bar K_{jl}}\\
 \label{stima01}
 & \leq & \displaystyle{-2\sum_{j,l}\left(\lambda_j-\lambda_l\right)^2}=-4m\aao.
\end{eqnarray}
Since $2/a \geq 2m-2 \geq m+3 = \rr$, we have
$$
2b\left(\aaa+\rr\right)-4m\left(\aaa-\frac{1}{m}\hhh\right)  =  -\frac 4a (\aaa - a \hhh - \frac a2 b\rr)  \leq  -\frac 4a Q.
$$
By the maximum principle, the assertion follows.
\cvd

\begin{Proposition}\label{pinching_preserved_codim}
Let $\mm_0$ be a closed submanifold of $\cc\pp^n$ of dimension $m$ and codimension $\displaystyle{2\leq k<\frac{2n-3}{5}}$. Then the pinching condition
\begin{equation}\label{pinching_stretto_codim}
\aaa \leq \frac{1}{m-1+\ep} \hhh+\frac{m-3-4k}{m}(1-\ep)
\end{equation}
is preserved by the flow for  any $\ep \in [0,1).$
\end{Proposition}
\proof 
Again, let us set
$Q=\aaa - a  \hhh -b$, where
$$a=\frac{1}{m-1+\ep}, \qquad b= \frac{m-3-4k}{m}(1-\ep).
$$
 By Lemma \ref{evoluzione_AH} we have
\begin{equation}\label{evoluzione_Q}
\dt Q = \Delta Q-2(\left|\nb A\right|^2-a\left|\nb H\right|^2)+2R_1-2aR_2+P_a,
\end{equation}
where $P_{a}=I+II+III$, with
$$
\begin{array}{c}
I=\displaystyle{4\sum_{i,j,p,q}\bar R_{ipjq}\left(\sum_{\alpha}h_{pq}^{\alpha}h_{ij}^{\alpha}\right)-4\sum_{j,s,p}\bar R_{sjsp}\left(\sum_{i,\alpha}h_{pi}^{\alpha}h_{ij}^{\alpha}\right)},\\
II=\displaystyle{2\sum_{s,\alpha,\beta}\bar R_{s\alpha s\beta}\left(\sum_{ij}h_{ij}^{\alpha}h_{ij}^{\beta}\right)-2a\sum_{s,\alpha,\beta}\bar R_{s\alpha s\beta}H^{\alpha}H^{\beta}},\\
III=\displaystyle{-8\sum_{j,p,\alpha,\beta}\bar R_{jp\alpha\beta}\left(\sum_i h_{ip}^{\alpha}h_{ij}^{\beta}\right)}.
\end{array}
$$
By Lemma \ref{grad05} the gradient terms in equation \eqref{evoluzione_Q} are non-positive and it suffices to consider the contribution of the reaction terms. Let us divide the analysis into two cases: $H=0$ and $H\neq 0$. Consider first a point where $Q=0$ and $H \neq 0$. To estimate {\em I}, we fix $\alpha$ and choose a tangent basis
$\{\widetilde{e}_1, \dots ,\widetilde{e}_m\}$, not necessarily of kind (B1) or (B2), that diagonalizes $h^{\alpha}$, i.e. $h_{ij}^{\alpha}=\lambda_i^{\alpha}\delta_{ij}$. Like in estimate \eqref{stima01}, we have
$$
\begin{array}{l}
\displaystyle{4\sum_{i,j,p,q}\bar R_{ipjq}h_{pq}^{\alpha}h_{ij}^{\alpha}-4\sum_{j,s,p}\bar R_{sjsp}\left(\sum_ih_{pi}^{\alpha}h_{ij}^{\alpha}\right)}\\
=4\displaystyle{\sum_{i,p}\bar R_{ipip}}(\lambda_i^{\alpha}\lambda_p^{\alpha}-(\lambda_i^{\alpha})^2)\\
=-2\displaystyle\sum_{i,p}\bar K_{ip}(\lambda_i^{\alpha}-\lambda_p^{\alpha})^2
\leq-4m\accentset{\circ\phantom{^{\alpha}}}{\left|h^{\alpha}\right|}^2.
\end{array}
$$
Hence we obtain
\begin{equation}\label{I}
I\leq-4m\aao.
\end{equation}

A basis of type (B2) is useful for estimating the terms $II$ and $III$. We recall that the curvature tensor of the Fubini-Study metric, for every $X$, $Y$, $Z$ and $W$ tangent vector fields of $\cc\pp^n$, is
\begin{eqnarray}\label{barR}
\nonumber\bar R(X,Y,Z,W) & = & \gfs(X,Z)\gfs(Y,W)-\gfs(X,W)\gfs(Y,Z)\\
& &+ \gfs(X,JZ)\gfs(Y,JW)-\gfs(X,JW)\gfs(Y,JZ)\\
\nonumber& &+2\gfs(X,JY)\gfs(Z,JW).
\end{eqnarray}
In order to study the term $II$, note that, with our choice of the basis, we have that $\bar R_{s\alpha s\beta}=0$ for any $s$ if $\alpha\neq\beta$. Otherwise we have
$$
\bar R_{s\alpha s\alpha}=\bar K_{s\alpha}=1+3\gfs(e_s,Je_{\alpha})^2,
$$
which implies that $1 \leq \bar K_{s\alpha} \leq 1+3\delta_{s,\alpha-m}$. Therefore, since $a\geq \frac 1m$, we have
\begin{eqnarray}\label{II}
\nonumber II & = &2\sum_{s,\alpha}\bar K_{s\alpha}\left(\left|h^\alpha\right|^2-a\left|H^\alpha\right|^2\right)\\
 \nonumber & = & 2\sum_{s,\alpha}\bar K_{s\alpha} \left(\left|\accentset{\circ\ }{h}^\alpha\right|^2-\left(a-\frac 1m\right)\left|H^\alpha\right|^2\right)\\
\nonumber & \leq & 2\sum_{s,\alpha}(1+3\delta_{s,\alpha-m})\left|\accentset{\circ\ }{h}^\alpha\right|^2 
\\
 & = & 2(m+3)\aao.
\end{eqnarray}

The most difficult term is $III$. 
Since $\bar R_{jp\alpha\beta}$ is anti-symmetric in $j,p$, while $h_{jp}^{\alpha}$ is symmetric, we have
\begin{eqnarray*}
III & =& \displaystyle{-8\sum_{j,p,\alpha,\beta}\bar R_{jp\alpha\beta}\left(\sum_i h_{ip}^{\alpha}h_{ij}^{\beta}\right)}
=\displaystyle{-8\sum_{j,p,\alpha,\beta}\bar R_{jp\alpha\beta}\left(\sum_i \accentset{\circ}{h}_{ip}^{\alpha}\accentset{\circ}{h}_{ij}^{\beta}\right)}.
\end{eqnarray*}
We now analyze the possible values of $\bar R_{jp\alpha\beta}$. First fix $\alpha$ and $\beta$ coupled by \eqref{base01}, meaning that $\min\{\alpha,\beta\}=m+2r-1$ and $\max\{\alpha,\beta\}=m+2r$ for some $r \leq k/2$. By symmetry, it suffices to consider the case where $\alpha<\beta$. We find 
\begin{eqnarray*}
\bar R_{jp\alpha\beta} & = & \tau_r^2\left(\delta_{j,2r-1}\delta_{p,2r}-\delta_{j,2r}\delta_{p,2r-1}\right) -2\nu_r\gfs(e_j,Je_p),
\end{eqnarray*}
and
$$
\gfs(e_j,Je_p)=\left\{\begin{array}{llll}
-\nu_s & \text{if } j=2s, & p=2s-1, &s\leq\frac k2; \medskip \\
\nu_s & \text{if } j=2s-1, & p=2s, & s\leq\frac k2; \medskip\\
1 & \text{if } j=k+2s, & p=k+2s-1, \quad &s\leq \frac {m-k}{2}; \medskip\\
-1 & \text{if } j=k+2s-1, \quad & p=k+2s, & s\leq \frac {m-k}{2}; \medskip\\
0 & \text{otherwise.}
\end{array}
\right.
$$
If $\alpha$ and $\beta$ are not coupled by \eqref{base01}, there are two indices $r\neq s$ such that $\alpha$ is (or is coupled with) $e_{m+2r-1}$ and $\beta$ is (or is coupled with) $e_{m+2s-1}$. In this case we have
$$
\bar R_{jp\alpha\beta}=\tau_r\tau_s\left(\delta_{j,\alpha-m}\delta_{p,\beta-m}-\delta_{j,\beta-m}\delta_{p,\alpha-m}\right).
$$

Using what we have just found and summing all similar terms we have
\begin{eqnarray*}
III & = & 16\sum_r\left(2\nu_r^2-\tau_r^2\right)\sum_i\left(\accentset{\circ}{h}_{i\ 2r}^{m+2r-1}\accentset{\circ}{h}_{i\ 2r-1}^{m+2r}-\accentset{\circ}{h}_{i\ 2r-1}^{m+2r-1}\accentset{\circ}{h}_{i\ 2r}^{m+2r}\right)\\
 & & -8\sum_{r\neq s\leq\frac k2}\tau_r\tau_s\sum_i\left(\accentset{\circ}{h}_{i\ 2s}^{m+2r}\accentset{\circ}{h}_{i\ 2r}^{m+2s}-\accentset{\circ}{h}_{i\ 2r}^{m+2r}\accentset{\circ}{h}_{i\ 2s}^{m+2s}\right)\\
 & & -16\sum_{r\neq s, r \leq\frac k2, s\leq\frac {k+1}2}\tau_r\tau_s\sum_i\left(\accentset{\circ}{h}_{i\ 2s-1}^{m+2r}\accentset{\circ}{h}_{i\ 2r}^{m+2s-1}-\accentset{\circ}{h}_{i\ 2r}^{m+2r}\accentset{\circ}{h}_{i\ 2s-1}^{m+2s-1}\right)\\
 & & -8\sum_{r\neq s\leq\frac {k+1}2}\tau_r\tau_s\sum_i\left(\accentset{\circ}{h}_{i\ 2s-1}^{m+2r-1}\accentset{\circ}{h}_{i\ 2r-1}^{m+2s-1}-\accentset{\circ}{h}_{i\ 2r-1}^{m+2r-1}\accentset{\circ}{h}_{i\ 2s-1}^{m+2s-1}\right)\\
 & & +32\sum_{r\neq s\leq\frac k2}\nu_r\nu_s\sum_i\left(\accentset{\circ}{h}_{i\ 2s}^{m+2r-1}\accentset{\circ}{h}_{i\ 2s-1}^{m+2r}-\accentset{\circ}{h}_{i\ 2s-1}^{m+2r-1}\accentset{\circ}{h}_{i\ 2s}^{m+2r}\right)\\
& & +32\sum_{r \leq \frac k2}\nu_r\sum_{s\leq\frac {m-k}{2}}\sum_i\left(\accentset{\circ}{h}_{i\ k+2s-1}^{m+2r-1}\accentset{\circ}{h}_{i\ k+2s}^{m+2r}-\accentset{\circ}{h}_{i\ k+2s}^{m+2r-1}\accentset{\circ}{h}_{i\ k+2s-1}^{m+2r}\right).
\end{eqnarray*}
Obviously $III\leq\left|III\right|$. Using repeatedly the triangle inequality and Young's inequality, and taking into account that for any $r$ and $s$
$$
\left\{
\begin{array}{l}
\left|2\nu_r^2-\tau_r^2\right| \leq 2,\\
\left|\tau_r\tau_s\right|\leq 1,\\
\left|\nu_r\nu_s\right|\leq 1,\\
\left|\nu_r\right|\leq 1,
\end{array}\right.
$$
we have
\begin{eqnarray*}
III & \leq & 16 \sum_{r \leq \frac k2} \sum_i \left(\left|\accentset{\circ}{h}_{i\ 2r}^{m+2r-1}\right|^2+\left|\accentset{\circ}{h}_{i\ 2r-1}^{m+2r}\right|^2+\left|\accentset{\circ}{h}_{i\ 2r-1}^{m+2r-1}\right|^2+\left|\accentset{\circ}{h}_{i\ 2r}^{m+2r}\right|^2\right)\\
 & & +4\sum_{r\neq s\leq\frac k2} \sum_i \left(\left|\accentset{\circ}{h}_{i\ 2s}^{m+2r}\right|^2+\left|\accentset{\circ}{h}_{i\ 2r}^{m+2s}\right|^2+\left|\accentset{\circ}{h}_{i\ 2r}^{m+2r}\right|^2+\left|\accentset{\circ}{h}_{i\ 2s}^{m+2s}\right|^2\right)\\
 & & +8\sum_{r\neq s, r \leq\frac k2, s\leq\frac {k+1}2} \sum_i \left(\left|\accentset{\circ}{h}_{i\ 2s-1}^{m+2r}\right|^2+\left|\accentset{\circ}{h}_{i\ 2r}^{m+2s-1}\right|^2+\left|\accentset{\circ}{h}_{i\ 2r}^{m+2r}\right|^2+\left|\accentset{\circ}{h}_{i\ 2s-1}^{m+2s-1}\right|^2\right)\\
 & & +4\sum_{r\neq s\leq\frac {k+1}2} \sum_i \left(\left|\accentset{\circ}{h}_{i\ 2s-1}^{m+2r-1}\right|^2+\left|\accentset{\circ}{h}_{i\ 2r-1}^{m+2s-1}\right|^2+\left|\accentset{\circ}{h}_{i\ 2r-1}^{m+2r-1}\right|^2+\left|\accentset{\circ}{h}_{i\ 2s-1}^{m+2s-1}\right|^2\right)\\
& & +16\sum_{r\neq s\leq\frac k2} \sum_i \left(\left|\accentset{\circ}{h}_{i\ 2s}^{m+2r-1}\right|^2+\left|\accentset{\circ}{h}_{i\ 2s-1}^{m+2r}\right|^2+\left|\accentset{\circ}{h}_{i\ 2s-1}^{m+2r-1}\right|^2+\left|\accentset{\circ}{h}_{i\ 2s}^{m+2r}\right|^2\right)\\
& & +16\sum_{r \leq \frac k2,s\leq\frac {m-k}{2}} \sum_i \left(\left|\accentset{\circ}{h}_{i\ k+2s-1}^{m+2r-1}\right|^2+\left|\accentset{\circ}{h}_{i\ k+2s}^{m+2r}\right|^2+\left|\accentset{\circ}{h}_{i\ k+2s}^{m+2r-1}\right|^2+\left|\accentset{\circ}{h}_{i\ k+2s-1}^{m+2r}\right|^2\right).
\end{eqnarray*}
Note that, if $k=2$, there are no indices $r\neq s \leq \frac{k+1}{2}$. Then, some of the sums in the expressions above are empty and we easily find that
$$
III\leq 16\aao.
$$
If $k>2$, by collecting similar terms we find
\begin{eqnarray*}
III & \leq & \sum_{i,r}\left(16\left|\accentset{\circ}{h}_{i\ 2r}^{m+2r-1}\right|^2+16\left|\accentset{\circ}{h}_{i\ 2r-1}^{m+2r}\right|^2+8k\left|\accentset{\circ}{h}_{i\ 2r}^{m+2r}\right|^2+8k\left|\accentset{\circ}{h}_{i\ 2r-1}^{m+2r-1}\right|^2\right)\\
 & & +24\sum_{i,r\neq s\leq\frac k2}\left(\left|\accentset{\circ}{h}_{i\ 2s}^{m+2r}\right|^2+\left|\accentset{\circ}{h}_{i\ 2s}^{m+2r-1}\right|^2+\left|\accentset{\circ}{h}_{i\ 2s-1}^{m+2r}\right|^2+\left|\accentset{\circ}{h}_{i\ 2s-1}^{m+2r-1}\right|^2\right)\\
& & +16\sum_{i,r,s\leq \frac {m-k}{2}}\left(\left|\accentset{\circ}{h}_{i\ k+2s-1}^{m+2r-1}\right|^2+\left|\accentset{\circ}{h}_{i\ k+2s}^{m+2r}\right|^2+\left|\accentset{\circ}{h}_{i\ k+2s}^{m+2r-1}\right|^2+\left|\accentset{\circ}{h}_{i\ k+2s-1}^{m+2r}\right|^2\right)\\
 & \leq & 8k\aao.
\end{eqnarray*}
So we can say that in any case 
\begin{equation}\label{III}
III\leq 8k\aao.
\end{equation}
By \eqref{I}, \eqref{II} and \eqref{III}, we conclude that
$$P_a=I+II+III\leq -2(m-3-4k)\aao.$$ 

Now let $R=2R_1-2aR_2+P_a$. If we again consider a frame of type (B1), Lemma \ref{l.est} says that at any point with $Q=0$ we have
\begin{eqnarray*}
R & \leq &\left(6 - \frac{2}{ma-1}  \right) \aao\amo  + \left(\frac{2mab}{ma-1}-2(m-3-4k)\right) \auo -3\amoq   \\
& &+ \left(\frac{4b}{ma-1} -2(m-3-4k)\right)\amo- \frac{2b^2}{ma-1}. \\
\end{eqnarray*}
Observe that, for our choice of $a$ and $b$, the coefficient of $\aao\amo$ is negative, while the one multiplying $\auo$ is zero. 
In addition, the assumptions $Q=0$ and $a>1/m$ imply that $\aao\geq b$.
Using this, we obtain
\begin{eqnarray*}
R &   \leq &
-3\amoq   + \left[\left(6 - \frac{2}{ma-1}\right)b +\frac{4b}{ma-1} -2(m-3-4k)\right]\amo- \frac{2b^2}{ma-1} \\
& = & -3\amoq+4b\amo+2b(b-m+3+4k).
\end{eqnarray*}
Using $4b\amo \leq 3\amoq+\frac43b^2$, we deduce
$$
R\leq 2b \left(\frac 53 b-m+3+4k \right).
$$
Our choice of $b$ then implies that $R < 0$.
 
Finally, let us consider the case of a point where $Q=\hhh=0$. Then we have $\aaa=\aao=b$, $R_2=0$. Moreover, using Theorem 1 from \cite{LL}, we find that $2R_1\leq 3\aaaq=3b^2$. As before, we obtain that $P_a\leq -2(m-3-4k)\aao=-2(m-3-4k)b$. Therefore,
$$
R\leq 3b^2-2(m-3-4k)b,
$$
which is negative for our choice of  $b$. By the maximum principle, the assertion follows. \cvd

\section{The traceless second fundamental form}\setcounter{equation}{0}\setcounter{Theorem}{0}
Following an approach which goes back to \cite{Ha,H1}, the description of the asymptotic behavior of $\mm_t$ will be obtained analyzing the traceless part of the second fundamental form and showing that it becomes small in a suitable sense if the curvature becomes unbounded. 

Since our initial manifold $\mm_0$ satisfies the assumption \eqref{pinching_codim}, it also satisfies inequality \eqref{pinching_stretto_hyp}, respectively \eqref{pinching_stretto_codim}, for some $\ep>0$. We know from the results of the previous section that these inequalitis are preserved by the flow for all $t>0$. 

As in \cite{H3,Ba}, we introduce the functions 
$$W:=\alpha\hhh+\beta, \qquad \fs:=\frac{\aao}{W^{1-\sigma}}.$$
Here $\sigma$ is a suitably small non-negative constant, while $\alpha,\beta$ are defined by
\begin{equation}\label{alphabeta}
\alpha=\left\{ \begin{array}{ll}
\displaystyle \frac{2}{(m-1 +\ep)(2+\rr -2\ep)} &\mbox{ if }k=1 \medskip \\
\displaystyle \frac{m-10}{3m^2} & \mbox{ if }k \geq 2
\end{array}
\right. \qquad
\beta=\left\{ \begin{array}{ll}
2 &\mbox{ if }k=1 \medskip \\
\displaystyle \frac{m-3-4k}m & \mbox{ if }k \geq 2.
\end{array}
\right.
\end{equation}

The main result of this section is the next proposition, which gives a differential inequality satisfied by $f_{\sigma}$.

\begin{Proposition}\label{fs_codim}
Under the assumptions of Theorem
\ref{maincodim} there is a $\sigma_1$ depending only on $\mm_0$ that for all $0\leq\sigma\leq\sigma_1$
\begin{equation}\label{dt_fs_codim}
\dt f_{\sigma}\leq \Delta f_{\sigma}+\frac{2\alpha (1-\sigma)}{W}\la\nb f_{\sigma},\nb \left|H\right|^2\ra-2C_1W^{\sigma-1}\left|\nb H\right|^2+2\sigma\aaa f_{\sigma}-2C_2f_{\sigma},
\end{equation}
for some constants $C_1>0$ and $C_2>0$ depending only on $m$ and the initial data.
\end{Proposition}
\proof 
Let us analyze the evolution equation for $f_\sigma$.
A straightforward computation gives
\begin{equation}\label{delta_fs}
\begin{array}{rcl}
\Delta \fs & = &\displaystyle{W^{\sigma-1}\Delta\aao-\alpha(1-\sigma)\frac{\fs}{W}\Delta\hhh-\frac{2\alpha(1-\sigma)}{W}\la\nb\fs,\nb\hhh\ra}\\
 & & \displaystyle{+\alpha^2\sigma(1-\sigma)\frac{\fs}{W^2}\left|\nb\hhh\right|^2.}
\end{array}
\end{equation}
Therefore
\begin{eqnarray*}
\dt \fs - \Delta \fs & = & W^{\sigma -1}\left(\dt\aao -\Delta \aao\right) -\alpha(1-\sigma)\frac{\fs}{W}\left(\dt\hhh- \Delta \hhh\right) \\
& & +\frac{2\alpha(1-\sigma)}{W}\la\nb\fs,\nb\hhh\ra -\alpha^2\sigma(1-\sigma)\frac{\fs}{W^2}\left|\nb\hhh\right|^2.
\end{eqnarray*}

Let us first consider the case of hypersurfaces $k=1$. Using Lemma \ref{evoluzione_AH_hyp}, and neglecting the negative 
$\left|\nb\hhh\right|^2$ term, we have
\begin{equation}\label{eq_fs_hyp}
\begin{array}{rcl}
\dt\fs & \leq & \displaystyle{\Delta\fs+\frac{2\alpha(1-\sigma)}{W}\la\nb\fs,\nb\hhh\ra-2W^{\sigma-1}\left|\nb A\right|^2}\\
 & &\displaystyle{+2W^{\sigma-1}\left[\frac{1}{m}+f_0(1-\sigma)\alpha\right]\left|\nb H\right|^2}\\
 & &\displaystyle{ +2\beta\frac{(1-\sigma)}{W}\fs\left(\aaa+\rr\right)+2\sigma\fs\left(\aaa+\rr\right)}\\
 & & \displaystyle{-4W^{\sigma-1}\left(h_{ij}h_j^{\phantom j p}\bar R_{pli}^{\phantom{mli}l}-h^{ij}h^{lp}\bar R_{pilj}\right).}
\end{array}
\end{equation}

Our choice of $\alpha$ and $\beta$ gives $0\leq f_0<1$. Hence, by Lemma \ref{grad05},
\begin{equation}\label{eq_pe_01}
\begin{array}{rcl}
-\left|\nb A\right|^2 &+&\displaystyle{\left[\frac{1}{m}+f_0(1-\sigma)\alpha\right]\left|\nb H\right|^2}\\
 &\leq&\displaystyle{ \left(\frac{1}{m}+\alpha\right)\left|\nb H\right|^2-\left|\nb A\right|^2
 = -C_1\left|\nb H\right|^2,}
\end{array}
\end{equation}
where $C_1=\frac{3}{m+2}-\frac 1m-\alpha$ is positive for our choice of $\alpha$ and $m \geq 5$.
It remains to estimate the reaction terms. Let us set
\begin{eqnarray*}
R & := &2\beta\frac{(1-\sigma)}{W}\fs\left(\aaa+\rr\right)+2\sigma\fs\left(\aaa+\rr\right) -4W^{\sigma-1}\left(h_{ij}h_j^{\phantom j p}\bar R_{pli}^{\phantom{mli}l}-h^{ij}h^{lp}\bar R_{pilj}\right).
\end{eqnarray*}
Using inequality (\ref{stima01}) we have
$$
R\leq 2\fs\left[\beta(1-\sigma)\frac{\aaa+\rr}{W}+\sigma(\aaa+\rr)-2m\right].
$$
From  \eqref{pinching_stretto_hyp} and the definitions   \eqref{barerre}, \eqref{alphabeta} of $\bar r$, $\alpha$ and $\beta$, we obtain
$$
\aaa+\rr\leq \frac{1}{m-1+\ep}\hhh+2(1-\ep)+\rr = \frac{2+\rr-2\ep}{\beta}W.
$$
Since $m \geq 5$ and $\ep$ is small, we have
\begin{eqnarray*}
R&\leq& 2 \fs \left[(1-\sigma)(\beta+\rr -2\ep)+\rr\sigma-2m\right] +2\sigma \fs \aaa\\
& = & 2\fs\left[5-m-2\ep+\sigma(2\ep-2)\right] +2\sigma \fs \aaa \leq-4\ep\fs + 2\sigma \fs \aaa.
\end{eqnarray*}
This inequality, together with \eqref{eq_fs_hyp} and \eqref{eq_pe_01}, implies the assertion for the case of hypersurfaces, with $C_2=2\ep$. 

Let us now turn to the case $k\geq 2$. From Lemma \ref{evoluzione_AH} and the properties of the curvature tensor $\bar R$, arguing as in the estimation of term $II$ in the proof of Proposition \ref{pinching_preserved_codim}, we find
\begin{eqnarray}
\nonumber\dt \hhh  &=&  \Delta\hhh -2\left|\nb H\right|^2+2R_2+2\sum_{s,\alpha} \bar K_{s\alpha}\left|H^\alpha\right|^2\\
 \label{H_bis}& \geq & \Delta\hhh -2\left|\nb H\right|^2+2R_2+2m\hhh.
\end{eqnarray}
Moreover, by Lemma \ref{evoluzione_AH}, we have
$$
\dt\aao=\Delta\aao-2\left(\nba-\frac 1m\nbh\right)+2\left(R_1-\frac 1m R_2\right)+P_{\frac 1m},
$$
where, like in the proof of Proposition \ref{pinching_preserved_codim},
$$
P_{\frac 1m}\leq -2(m-3-4k)\aao.
$$
Then
\begin{eqnarray*}
\dt\fs & \leq & W^{\sigma-1}\left(\Delta\aao-2\left(\nba-\frac 1m\nbh\right)\right)\\
 & & +W^{\sigma-1}\left(2\left(R_1-\frac 1m R_2\right)-2(m-3-4k)\aao\right)\\
 & &  -\alpha(1-\sigma)\frac{\fs}{W}\left(\Delta\hhh -2\left|\nb H\right|^2+2R_2+2m\hhh\right).
\end{eqnarray*}
 Using the expression found previously for $\Delta\fs$,  we obtain
\begin{equation}\label{eq_fs_codim}
\begin{array}{rcl}
\dt\fs & \leq & \displaystyle{\Delta\fs+\frac{2\alpha(1-\sigma)}{W}\la\nb\fs,\nb\hhh\ra-2W^{\sigma-1}\left|\nb A\right|^2}\\
 & &\displaystyle{+2W^{\sigma-1}\left[\frac{1}{m}+f_0(1-\sigma)\alpha\right]\left|\nb H\right|^2}\\
 & & \displaystyle{+2W^{\sigma-1}\left(R_1-\frac{1}{m}R_2\right)-2\alpha(1-\sigma)\frac{\fs}{W}R_2}\\
& & \displaystyle{-2m\alpha(1-\sigma)\frac{\fs}{W}\hhh-2(m-3-4k)W^{\sigma-1}}\aao.
\end{array}
\end{equation}
To estimate the gradient terms, we use Lemma \ref{grad05}. Let us set
$$
C_1= \frac{16}{9(m+2)} - \frac{4m-10}{3m^2},
$$
which is positive for all $m \geq 0$. Then we have, using again $0\leq f_0 <1$,
\begin{eqnarray*}
\left[\frac{1}{m}+f_0(1-\sigma)\alpha\right]\left|\nb H\right|^2 
 &\leq& \left(\frac{1}{m}+\alpha\right)\left|\nb H\right|^2 = \frac{4m-10}{3m^2} \nbh \\
&=&  \left( \frac{16}{9(m+2)} -C_1 \right) |\nb H|^2 \leq |\nb A |^2 - C_1\nbh,
\end{eqnarray*}
which yields the desired estimate. Let us now analyze the reaction terms. We can write them as
$$
R= 2W^{\sigma-2} R' + 2\alpha\sigma\frac{\fs}{W}R_2
$$
where
$$
R'=\left(R_1-\frac{1}{m}R_2\right)W-\alpha\aao R_2-\alpha m(1-\sigma)\aao\hhh-(m-3-4k)\aao W.
$$
We first estimate
\begin{equation}\label{eq_pe_02}
\begin{array}{rcl}
\displaystyle{2\alpha\sigma\frac{\fs}{W}R_2} & \leq &\displaystyle{2\alpha\sigma\frac{\fs}{W}\aaa\hhh}\\
& = &\displaystyle{2\sigma\fs\aaa-2\sigma\beta\frac{\fs}{W}\hhh}\\
&\leq& \displaystyle{2\sigma\fs\aaa.}
\end{array}
\end{equation}
The rest of the proof is devoted to the estimation of $R'$. By Lemma \ref{l.est}
$$
R_1-\frac 1m R_2\leq\auoq+\frac{1}{m}\auo\hhh+4\auo\amo+\frac{3}{2}\amoq.
$$
Moreover $\aao=\auo+\amo$ and  $R_2=\auo\hhh+\frac{1}{m}\hhhq$, so
\begin{eqnarray}
\nonumber R' & \leq & 3\alpha\auo\amo\hhh+\frac 32\alpha\amoq\hhh-\frac{\alpha}{m}\amo\hhhq\\
 \nonumber& & +\beta\auoq+4\beta\auo\amo+\frac 32\beta\amoq\\
 \label{R'}& & +\left(\frac{\beta}{m}-m\alpha(1-\sigma)-\alpha(m-3-4k)\right)\auo\hhh\\
 \nonumber& & -\alpha\left(m(1-\sigma)+m-3-4k\right)\amo\hhh\\
 \nonumber& & -\beta(m-3-4k)\left(\auo+\amo\right).
\end{eqnarray}

\noindent Since the pinching condition \eqref{pinching_codim} holds, we have that 
$$
\left(\frac{1}{m-1}-\frac 1m\right)\hhh\geq\left(\auo+\amo-\beta \right).
$$
Then we have 
\begin{eqnarray*}
R' & = & R'+\frac{3\alpha}{m(m-1)}\amo\hhhq+\frac{2\beta}{m(m-1)} \left(\auo+\amo\right)\hhh\\
 & & -\frac{3 \alpha}{m(m-1)}
 \amo\hhhq-\frac{2\beta}{m(m-1)}\left(\auo+\amo\right)\hhh\\
 &\leq & R'+\frac{3 \alpha}{m(m-1)}\amo\hhhq+\frac{2\beta}{m(m-1)}\left(\auo+\amo\right)\hhh\\
 & & -3\alpha(\auo+\amo-\beta)\amo\hhh-2\beta\left(\auo+\amo\right)\left(\auo+\amo-\beta\right) \\
 & \leq & {\ro-\alpha\frac{m-4}{m(m-1)}} \amo \hhhq \\
 & & + \left[ \beta\left( \frac 2{m-1}- \frac 1m \right) -\alpha(m(2-\sigma)-3-4k) \right] \auo \hhh \\
 & & + \left[ \beta\left( \frac2{m-1} -  \frac 2m \right) -\alpha(m(2-\sigma)-3-4k-3\beta) \right] \amo \hhh \\
 & & +\beta\left(2 \beta-m+3+4k\right)\left(\auo+\amo\right).
\end{eqnarray*}
Our hypotheses give $\beta \leq\frac{1}{4}(m-3-4k)$. We can further assume that $\sigma$ is small, say $\sigma<\frac 14$. 
Using these inequalities, the condition $m>4k+3$ and the inequalities
$$\frac{2\beta}{m(m-1)} < \frac{\beta (m+1)}{m(m-1)} < \frac {m-3-4k}{4m},$$
we obtain
\begin{eqnarray*}
R'&\leq& \left[\beta\frac {m+1}{m(m-1)}-\frac{\alpha}{4}\left(7m-12-16k\right)\right]\auo\hhh\\
 & & +\left[ \frac{2\beta}{m(m-1)}  - \alpha \left( \frac 74 m -3-4k -3 \beta\right)\right]\amo\hhh\\
 & & +\beta\left(2\beta-m+3+4k\right)\left(\auo+\amo\right) \\
 & \leq & \left[\frac{m-3-4k}{4m}-\frac{\alpha}{4}\left(7m-12-16k\right)\right]\auo\hhh\\
 & & +\left[  \frac{m-3-4k}{4m}  -\frac \alpha 4 \left( 4 m -3-4k\right)\right]\amo\hhh\\
 && -2\beta^2 \left(\auo+\amo\right) \\
  & \leq & \left[\frac{m-11}{4m}-\frac{3m}{4}\alpha \right](\auo+\amo)\hhh\\
 && -2\beta^2 \left(\auo+\amo\right) \\
 &\leq&-C_2\aao W, 
\end{eqnarray*}
for some positive constant $C_2$  depending only on $m$. Together with \eqref{eq_pe_02}, this implies that
$$
R \leq 2\sigma\fs\aaa+2W^{\sigma-2}R'\leq 2\sigma\fs\aaa-2C_2\fs,
$$
which concludes our proof.
\cvd

We now prove some other estimates which will be needed in the following.

\begin{Lemma}\label{lemma_GE_codim}
We have the estimates:
\begin{enumerate}
\item [1)]$\dt \aao\leq \Delta\aao- 2C_3\nba+4\aaa\aao$, for some $C_3>0$ only depending on $m$,
\item [2)]$\dt\left| H\right|^4\geq\Delta\left| H\right|^4-12\hhh\nbh+\frac 4m\left| H\right|^6$.
\end{enumerate}
\end{Lemma}
\proof In the case of hypersurfaces, inequality 1) follows easily from Lemma \ref{evoluzione_AH_hyp}, inequality \eqref{grad02} and estimate \eqref{stima01}. For higher codimension we use Lemma \ref{evoluzione_AH}:
\begin{eqnarray*}
\dt\aao&\leq&\Delta\aao-2\left(\nba-\frac{1}{m}\nbh\right)\\
 & & +2\left(R_1-\frac{1}{m}R_2\right)+P_{\frac 1m}.
\end{eqnarray*}
By  Lemma \ref{grad05}, we have that $-2\left(\nba-\frac{1}{m}\nbh\right)\leq-2C_3\nba$ for some positive constant $C_3$.
Moreover, using Lemma \ref{l.est}, 
\begin{eqnarray*}
R_1-\frac{1}{m}R_2&\leq&\auoq+4\auo\amo+{\ro\frac 32}\amoq+\frac 1m \auo\hhh\\
& \leq& 2\left(\auo+\amo\right)^2+\frac 2m\hhh\left(\auo+\amo\right)=2\aao\aaa.
\end{eqnarray*}
Finally, like in the proof of Proposition \ref{pinching_preserved_codim},
$$
P_{\frac 1m}\leq-2(m-3-4k)\aao\leq0.
$$
This proves inequality 1). To prove the second part, we use again Lemma \ref{evoluzione_AH_hyp} and \eqref{evoluzione_AH}.
For hypersurfaces we obtain
\begin{eqnarray*}
\dt \hhhq & = & \Delta\left|H\right|^4-2\left|\nb\hhh\right|^2-4\hhh\nbh+4\hhhq(\aaa+\rr)\\
 & \geq & \Delta\hhhq-12\hhh\nbh+\frac 4m\hhhc.
\end{eqnarray*}
For higher codimension we use the inequality
$$
2R_2=2\hhh\left(\auo+\frac 1m\hhh\right)\geq\frac 2m \hhhq
$$
and we find
\begin{eqnarray*}
\dt\left|H\right|^4 & = & \Delta\left|H\right|^4-2\left|\nb\hhh\right|^2-4\hhh\nbh\\
 & & +2\hhh\left(2R_2+2\sum_{s,\alpha} \bar K_{s\alpha}\left|H^\alpha\right|^2\right)\\
 & \geq &  \Delta\left|H\right|^4-12\hhh\nbh+\frac{4}{m}\left|H\right|^6.
\end{eqnarray*}\cvd

Finally, we consider the evolution equation for $\nbh$. With the same proof of Corollary 5.10 in \cite{Ba}, we have the following result.

\begin{Proposition}\label{evoluzione_nbh}
There exists a constant $C_4$ depending only on $\mm_0$ such that
$$
\dt\nbh \leq \Delta\nbh+C_4(\hhh+1)\nba.
$$
\end{Proposition}

\section{Finite maximal time}\setcounter{equation}{0}\setcounter{Theorem}{0}

In this section we consider the case that our flow develops a singularity in finite time and prove convergence to a round point as stated in Theorem \ref{maincodim}.

Since $\mm_0$ is compact,  there is also an $\ep>0$ small enough such that 
\begin{equation}\label{pinching_codim_ep}
\aaa \leq a \hhh + b,
\end{equation}
where
\begin{equation}\label{pinching_codim_ep_1}
a=\frac{1}{m-1+\ep}, \qquad b=\left\{ \begin{array}{ll}2(1-\ep) \ & \mbox{ if }k=1 \medskip\\
\displaystyle\frac{m-3-4k}{m}(1-\ep) \ & \mbox{ if } k \geq 2.
\end{array}\right.
\end{equation}
We know that inequality \eqref{pinching_codim_ep} with the above choice of constants remains preserved during the flow. 
As in the previous section, we let $W=\alpha \hhh + \beta$, where $\alpha,\beta$ are chosen according to \eqref{alphabeta}. We observe that
\begin{eqnarray}
2mW & \geq & 2m \left( a-\frac 1m\right) \hhh + 2mb \nonumber \\
& = & 2 \frac{1-\ep}{m-1+\ep} \hhh + 2mb > a \hhh +  b \geq \aaa.
\label{amw}
\end{eqnarray}

\begin{Theorem}\label{PE_fin_codim}
Let the assumptions of Theorem \ref{maincodim} hold. If $\tm$ is finite, there are constants $C_0<\infty$ and $\sigma_0>0$ depending only on the initial manifold $\mm_0$ such that for all $0\leq t<\tm$ we have
$$
\aao\leq C_0(\hhh+1)^{1-\sigma_0}.
$$
\end{Theorem}
To prove this result we will bound from above the function $f_\sigma$ introduced in the previous section. For $\sigma > 0$, the  positive term $2\sigma\fs\aaa$ in (\ref{fs_codim}) prevents us from using the maximum principle. Therefore, as in Huisken \cite{H3} and Baker \cite{Ba}, we will obtain integral estimates on $f_\sigma$ exploiting the good negative $\left|\nb H\right|^2$ term by the divergence theorem. These estimates allow to deduce the desired sup-estimate through a standard iteration procedure. 

The starting point of our proof is the contracted Simons identity computed in \cite{AB}, formula (23). Using this we easily obtain
$$
\Delta\aao\geq 2 \accentset{\phantom{aaai}\circ}{\left|\nb{A}\right|}^2+2\la\accentset{\circ\phantom{ij}}{h_{ij}},\nb_i\nb_jH\ra+2Z-c|A|^2,
$$
where $c>0$ is a suitable constant only depending on $m,k$ and
\begin{equation*}
Z  =  \displaystyle{\sum_{i,j,p,\alpha,\beta}H^{\alpha}h_{ip}^{\alpha}h_{pj}^{\beta}h_{ij}^{\beta}-\sum_{\alpha,\beta}\left(\sum_{i,j}h_{ij}^{\alpha}h_{ij}^{\beta}\right)^2-\sum_{i,j,\alpha,\beta}\left(\sum_p\left(h_{ip}^{\alpha}h_{pj}^{\beta}-h_{jp}^{\alpha}h_{ip}^{\beta}\right)\right)^2.             }
\end{equation*} 
Using our pinching assumption we also deduce 
\begin{equation}\label{delta_A0}
\Delta\aao\geq 2 \accentset{\phantom{aaai}\circ}{\left|\nb{A}\right|}^2+2\la\accentset{\circ\phantom{ij}}{h_{ij}},\nb_i\nb_jH\ra+2Z-\gamma W,
\end{equation}
where $\gamma$ only depends on $m,k$.

To understand the properties of $Z$ in the case of hypersurfaces, it is interesting to relate the pinching condition \eqref{pinching_codim} to the positivity of the intrinsic sectional curvature of the submanifold $\mm_t$.  

\begin{Proposition}\label{sez_pos}
There exists a constant $c=c(m)$ such that if $k=1$ the intrinsic sectional curvature of $\mm_t$ satisfies at any point
$$
K>\ep c W>0.
$$
\end{Proposition}
\proof Let $e_1,\dots,e_m$ be a orthonormal tangent basis that diagonalizes the second fundamental form. For any $i\neq j$ the Gauss equation gives
$$
K_{ij}=\bar K_{ij}+\lambda_i\lambda_j.
$$
Like in \cite{H3},  we can use the following algebraic property:  for any $i\neq j$
\begin{eqnarray}
\nonumber \left|A\right|^2-\frac{1}{m-1}|H|^2 &=&-2\lambda_i\lambda_j+\left(\lambda_i+\lambda_j-\frac{|H|}{m-1}\right)^2+\sum_{l\neq i,j}\left(\lambda_l-\frac{|H|}{m-1}\right)^2\\
 \label{prop_alg}&\geq&-2\lambda_i\lambda_j.
\end{eqnarray}
Then we have
\begin{eqnarray*}
2K_{ij} & \geq & 2-\aaa+\frac{1}{m-1}\hhh\\
 & \geq & \left(\frac{1}{m-1}-a \right)\hhh+2-b \\
 & = & \ep \left(\frac{1}{(m-1)(m-1+\ep)}\hhh+2 \right) \\
 & \geq & \ep c \left(\alpha \hhh+\beta\right)>0,
\end{eqnarray*}
for a suitable $c=c(m)$. \cvd

We cannot use the same argument in higher codimension because we cannot diagonalize simultaneously the tensors $h^{\alpha}$, for $\alpha=m+1,\dots,2n$. However, as a consequence of our other estimates, we will prove at the end of this section that  also in this case the sectional curvature of the evolving submanifold becomes positive for time large enough.

\begin{Lemma}\label{lemma_Z}
There exists $\rho>0$ depending only on $m,k$ such that $Z$ satisfies
$$
Z+2m b\aao\geq\rho \ep \aao W.
$$
\end{Lemma}
\proof
Let us first consider the case of hypersurfaces. Choosing a basis that diagonalizes the second fundamental form, using Gauss equations, Proposition \ref{sez_pos}  and $\bar K\leq 4$, we have
\begin{eqnarray*}
Z & = &\left(\sum_i\lambda_i\right)\left(\sum_i\lambda_i^3\right)-\left(\sum_i\lambda_i^2\right)^2\\
 & = & \sum_{i<j}\lambda_i\lambda_j\left(\lambda_i-\lambda_j\right)^2\\
 & = & \sum_{i<j}K_{ij}\left(\lambda_i-\lambda_j\right)^2-\sum_{i<j}\bar K_{ij}\left(\lambda_i-\lambda_j\right)^2\\
 &\geq & \ep c(m) W\aao-4m\aao=\ep c(m) W\aao -2bm\aao.
\end{eqnarray*}
For $k\geq 2$ we need to distinguish the cases $H=0$ and $H\neq 0$. Let us examine first the case $H\neq 0$. We use an estimate proved by Andrews and Baker, see page 384 in \cite{AB}, which gives 
$$
Z\geq-\frac{m}{2}\auoq-\frac{3}{2}\amoq-\frac{m+2}{2}\auo\amo+\frac{1}{2(m-1)}\left(\auo+\amo\right)\hhh.
$$
Since (\ref{pinching_codim_ep}) and (\ref{pinching_codim_ep_1}) hold, we have $\displaystyle{\hhh\geq \frac{m(m-1+\ep)}{1-\ep}\left(\auo+\amo-b\right)}$. Then
\begin{eqnarray*}
Z &\geq&-\frac{m}{2}\auoq-\frac{3}{2}\amoq-\frac{m+2}{2}\auo\amo\\
 & &+\frac{m}{2(1-\ep)}\left(\auo+\amo\right)\left(\auo+\amo-b\right)\\
 & = & \frac{\ep m}{2(1-\ep)}\auoq+\frac{m-3+3\ep}{2(1-\ep)}\amoq\\
 & & +\frac{m-2+\ep(m+2)}{2(1-\ep)}\auo\amo-\frac{m}{2(1-\ep)}b\aao.
\end{eqnarray*}
We may assume that $\ep>0$ is small enough in order to have $2m>\frac{m}{2(1-\ep)}$. Then the above estimate shows that there exists $\rho_1=\rho_1(m)>0$ such that 
$$
Z+2mb\aao\geq\ep \rho_1\aaoq.
$$
On the other hand, using the definition of $Z$ and estimating various terms with Peter-Paul's inequality, we find
$$
Z\geq\rho_2\aao\hhh-\rho_3\aaoq,
$$
for $\rho_2$ and $\rho_3$ depending on $m$. Combining these two inequalities we obtain for any $0\leq c\leq 1$
$$
Z+2mb\aao\geq c\left(\rho_2\aao\hhh-\rho_3\aaoq+2mb\aao\right)+(1-c)\left(\ep\rho_1\aaoq\right).
$$
Choosing $\bar c=\frac{\ep\rho_1}{\ep\rho_1+\rho_3}$ we have
$$
Z+2mb\aao\geq \bar c\left( \rho_2\hhh+2mb\right)\aao.
$$
The assertion follows for $\rho$ small enough.

When $H=0$ we have $\aaa=\aao\leq b$ and $W=\beta=b$. Using Theorem 1 in \cite{LL} we find
$$
Z \geq-\frac{3}{2}\left|A\right|^4 \geq-\frac{3}{2}b\aao. 
$$
Hence we have
$$
Z+2mb\aao\geq \left(2m-\frac{3}{2}\right)b\aao =\left(2m-\frac{3}{2}\right)\aao W \geq \ep\rho\aao W,
$$
provided $\rho>0$ is small enough. \cvd

Next we derive a Poincar\'e-type inequality on $f_\sigma$.

\begin{Proposition}\label{integral_estimates}
There exists a constant $C_5$ depending only on $m,k$ and $\mm_0$ such that, for any $p\geq 2$, $0<\sigma<1/4$ and $\eta>0$, we have
\begin{eqnarray*}
\ep \rho\int_{\mm_t}\fs^p W d\mu&\leq&\left(\eta(p+1)+5\right)\int_{\mm_t} W^{\sigma-1}\fs^{p-1}\left|\nb H\right|^2d\mu+\frac{p+1}{\eta}\int_{\mm_t} \fs^{p-2}\left|\nb\fs\right|^2d\mu\\
 & &+4mb\int_{\mm_t} \fs^pd\mu+{\ro\frac 1p}C_5^p.
\end{eqnarray*}
\end{Proposition}
\proof Plugging equation (\ref{delta_A0}) into (\ref{delta_fs}), we find
\begin{eqnarray*}
\Delta\fs &\geq & 2W^{\sigma-1}\accentset{\phantom{aaai}\circ}{\left|\nb{A}\right|}^2+2W^{\sigma-1}\la\accentset{\circ\phantom{{ij}}}{h_{ij}},\nb_i\nb_jH\ra+2W^{\sigma-1}Z-\gamma W^{\sigma} -\alpha(1-\sigma)\frac{\fs}{W}\Delta\hhh\\
 & &-\frac{2\alpha(1-\sigma)}{W}\la\nb_i\fs,\nb_i\hhh\ra+\alpha^2\sigma(1-\sigma)\frac{\fs}{W^2}\left|\nb\hhh\right|^2.
\end{eqnarray*}
The terms $2W^{\sigma-1}\accentset{\phantom{aaai}\circ}{\left|\nb{A}\right|}^2$ and $\alpha^2\sigma(1-\sigma)\frac{\fs}{W^2}\left|\nb\hhh\right|^2$ are positive, so we can omit them. Thanks to Lemma \ref{lemma_Z}, we have
\begin{eqnarray*}
\Delta\fs &\geq & 2W^{\sigma-1}\la\accentset{\circ\phantom{ij}}{h_{ij}},\nb_i\nb_jH\ra -\alpha(1-\sigma)\frac{\fs}{W}\Delta\hhh-\frac{2\alpha(1-\sigma)}{W}\la\nb \fs,\nb\hhh\ra\\
 & & +2 \ep\rho W^{\sigma}\aao-4mb\fs-\gamma W^{\sigma}.
\end{eqnarray*}

We multiply the above inequality by $\fs^{p-1}$ and integrate on $\mm_t$. All terms, except the last two negative ones, can be estimated as in Lemma 2.4 in \cite{H3} and Proposition 5.5 in \cite{Ba}. In this way we obtain, for any $\eta>0$,
\begin{eqnarray*}
2\ep\rho\int_{\mm_t}\fs^p Wd\mu & \leq & \left(\eta(p+1)+5\right)\int_{\mm_t} W^{\sigma-1}\fs^{p-1}\left|\nb H\right|^2d\mu+\frac{p+1}{\eta}\int_{\mm_t}\fs^{p-2}\left|\nb \fs\right|^2d\mu\\
& & +4mb\int_{\mm_t} \fs^{p}d\mu+\gamma\int_{\mm_t} W^{\sigma}\fs^{p-1}d\mu.
\end{eqnarray*}
In order to estimate the last term we use Young's inequality:
$$
\gamma W^{\sigma}\fs^{p-1}\leq\gamma W\left( \frac{r^p}{p}W^{(\sigma-1)p}+\frac{p-1}{p}r^{-\frac{p}{p-1}}\fs^p\right),\quad\forall r>0.
$$
Choose $r$ such that $\frac{p-1}{p}\gamma r^{-\frac{p}{p-1}}=\ep \rho$. Observe that $r$ is uniformly bounded from above for large $p$. Moreover $(\sigma-1)p+1<0$ and $W\geq\beta>0$. Then $W^{(\sigma-1)p+1}\leq\beta^{(\sigma-1)p+1}$ and we have

\begin{eqnarray*}
\frac{1}{p}\gamma r^p\int_{\mm_t} W^{(\sigma-1)p+1}d\mu & \leq &\frac{1}{p} \gamma r^p\beta^{(\sigma-1)p+1}{\rm vol}(\mm_t)\\
 & \leq & \frac{1}{p}\gamma r^p\beta^{(\sigma-1)p+1}{\rm vol}(\mm_0) \leq \frac{1}{p}C_5^p
\end{eqnarray*}
for a suitable $C_5>0$ depending on $\mm_0$. \cvd

We can now bound high $L^p$-norms of $\fs$, provided $\sigma$ is of order $p^{-\frac{1}{2}}$. This is the step where the hypothesis $\tm<\infty$ is used in an essential way.

\begin{Proposition}\label{Lp}
If $\tm<+\infty$, there is a constant $C_6$ depending only on $m,k,\mm_0,\tm$ such that for all 
$$
p\geq{\ro\frac{16}{C_1\ep}+2}\qquad \sigma\leq\frac{\sqrt{C_1}\rho}{2^7m\sqrt{p}}{\ro\ep^2}
$$
we have the inequality
$$
\left(\int_{\mm_t}\fs^pd\mu\right)^{\frac{1}{p}}\leq C_6, \qquad \mbox{ for all }t <\tm.
$$
\end{Proposition}
\proof
We multiply inequality (\ref{dt_fs_codim}) by $p\fs^{p-1}$, integrate and obtain
\begin{eqnarray*}
\frac{d}{dt}\int_{\mm_t}\fs^pd\mu & + &p(p-1)\int_{\mm_t}\fs^{p-2}\left|\nb\fs\right|^2d\mu+2C_1p\int_{\mm_t}\left|\nb H\right|^2W^{\sigma-1}\fs^{p-1}d\mu\\
& \leq & 4p\alpha\int_{\mm_t}\left|H\right|W^{-1}\left|\nb H\right|\left|\nb \fs\right|\fs^{p-1}d\mu+2\sigma p\int_{\mm_t}\aaa\fs^pd\mu\\
 & & -2C_2p\int_{\mm_t}\fs^pd\mu.
\end{eqnarray*}
Using that $\alpha\left|H\right|\leq W^{\frac{1}{2}}$ and $\fs\leq W^{\sigma}$, we have
\begin{eqnarray*}
\lefteqn{ 4p\alpha \int_{\mm_t}\left|H\right|W^{-1}\left|\nb H\right|\left|\nb \fs\right|\fs^{p-1}d\mu} \\
 &\leq& \frac{p(p-1)}{2}\int_{\mm_t}\left(\frac{\alpha\left|H\right|}{W}\fs\right)\fs^{p-2}\left|\nb\fs\right|^2W^{\frac{1}{2}-\sigma}d\mu\\
 & & +\frac{8p}{p-1}\int_{\mm_t}\frac{\alpha\left|H\right|}{W}W^{\sigma-\frac{1}{2}}\fs^{p-1}\left|\nb H\right|^2d\mu\\
&\leq &\frac{p(p-1)}{2}\int_{\mm_t}\fs^{p-2}\left|\nb\fs\right|^2d\mu\\
& & +\frac{8p}{p-1}\int_{\mm_t} W^{\sigma-1}\fs^{p-1}\nbh d\mu.
\end{eqnarray*}
With our choice of $p$, we have $C_1p\leq2C_1p-\frac{8p}{p-1}$. In addition, \eqref{amw} shows that $\aaa\leq 2mW$. Therefore
\begin{eqnarray*}
\frac{d}{dt}\int_{\mm_t}\fs^pd\mu & + &\frac{p(p-1)}{2}\int_{\mm_t}\fs^{p-2}\left|\nb\fs\right|^2d\mu+C_1p\int_{\mm_t}\left|\nb H\right|^2W^{\sigma-1}\fs^{p-1}d\mu\\
& \leq & 2\sigma p\int_{\mm_t}\aaa\fs^pd\mu -2C_2p\int_{\mm_t}\fs^pd\mu\\
& \leq & 4\sigma pm\int_{\mm_t} W\fs^pd\mu -2C_2p\int_{\mm_t}\fs^pd\mu.
\end{eqnarray*}

Thanks to Lemma \ref{integral_estimates}, we obtain for any $\eta>0$
\begin{eqnarray*}
\frac{d}{dt}\int_{\mm_t}\fs^pd\mu & + &\frac{p(p-1)}{2}\int_{\mm_t}\fs^{p-2}\left|\nb\fs\right|^2d\mu+C_1p\int_{\mm_t}\left|\nb H\right|^2W^{\sigma-1}\fs^{p-1}d\mu\\
& \leq & \frac{{\ro 4}\sigma pm}{{\ro \ep}\rho}\left[\left(\eta(p+1)+5\right)\int_{\mm_t} W^{\sigma-1}\fs^{p-1}\left|\nb H\right|^2d\mu\right.\\
 & &\phantom{AAAA}+\left.\frac{p+1}{\eta}\int_{\mm_t} \fs^{p-2}\left|\nb \fs\right|^2d\mu+4mb\int_{\mm_t} \fs^pd\mu+\frac 1pC^p_5\right]\\
 & & -2C_2p\int_{\mm_t}\fs^pd\mu.
\end{eqnarray*}
Choosing $\eta=\frac{\sqrt{C_1}{\ro\ep}}{4\sqrt{p}}$ and using our assumptions on $m$, $p$ and $\sigma$, we have
$$
 \frac{{\ro 4}\sigma pm}{{\ro\ep}\rho}\left(\eta(p+1)+5\right)\leq C_1p,\qquad \frac{{\ro 4}\sigma p(p+1)}{{\ro\ep}\rho\eta}\leq\frac{p(p-1)}{2}.
$$
Then
$$
\frac{d}{dt}\int_{\mm_t}\fs^pd\mu\leq \bar{C_2}\int_{\mm_t}\fs^pd\mu+ \bar{C_5},
$$
where
$$\displaystyle{\bar{C_2}=\frac{32m^2pb\sigma}{\rho}-2C_2p}, \qquad \displaystyle{\bar{C_5}=\frac{8\sigma m}{\rho}C^p_5}.$$
Since $\tm$ is finite, we obtain the assertion for a constant $C_6$ independent of $p$. \cvd

To prove Theorem \ref{PE_fin_codim}, we can now proceed as in \cite{H3} via a Stampacchia iteration procedure to uniformly bound the function $\fs$ when $\tm<\infty$.

Next we establish a gradient estimate for the mean curvature flow. This estimate is required to compare the mean curvature at different points of the submanifold. First we need some technical inequalities. As before, we denote by $C_i$ constants only depending on $m,k$ and the initial data.

\begin{Lemma}\label{H4fin}
\begin{eqnarray*}
\dt\hhh\aao &\leq&  \Delta(\hhh\aao)- C_3\hhh\nba+C_7\nba\\
&&+2\hhh\aao({\ro 3}\aaa+{\ro 4}m)
\end{eqnarray*}
for some constant $C_7>0$.
\end{Lemma}
  
\proof By Lemma \ref{evoluzione_AH} and \ref{lemma_GE_codim},
\begin{eqnarray*}
\dt\hhh\aao&\leq&\Delta(\hhh\aao)-2\la\nb\hhh,\nb\aao\ra - 2C_3\hhh\nba\\
 & &-{\ro 2}\aao\nbh+2\aao\hhh({\ro 3}\aaa+{\ro 4}m).
\end{eqnarray*}

Furthermore we have
\begin{eqnarray}
\nonumber-2\la\nb\hhh,\nb\aao\ra&\leq &4\left|H\right|\la\left|\nb H\right|,\nb\aao\ra\\
 \nonumber&\leq &8\left|H\right|\left|\nb H\right|\aao\left|\nb A\right|\\
 \nonumber& \leq & {\ro 6}\left|H\right|\sqrt{m+2}\nba\aao.
\end{eqnarray}

We can estimate the last term using Theorem \ref{PE_fin_codim} and Young inequality, to find that there exists a constant $C_7>0$ such that
\begin{eqnarray*}
{\ro 6}\left|H\right|\sqrt{m+2}\nba\aao &\leq &{\ro 6}\left|H\right|\sqrt{m+2}\nba\sqrt{C_0}\left(\hhh+1\right)^{\frac{1-\sigma}{2}}\\
 & \leq &  C_3\hhh\nba+C_7\nba.
\end{eqnarray*}
\cvd

\noindent Now we consider the function 
\begin{equation}\label{defg}
g=\hhh\aao+{\ro \left(\frac{C_7}{C_3}+1\right)}\aao.
\end{equation}
Using Lemma \ref{lemma_GE_codim}, Lemma \ref{H4fin} and $\hhh\leq m\aaa$ we obtain
\begin{eqnarray}
\nonumber\dt g & \leq & \Delta g- C_3\hhh\nba+C_7\nba+2\aao\hhh({\ro 3}\aaa+{\ro 4}m)\\
 \nonumber& & {\ro \left(\frac{C_7}{C_3}+1\right)}\left(- 2C_3\nba+4\aaa\aao\right)\\
 \nonumber&\leq & \Delta g-C_3\hhh\nba-4C_3\nba+2\aao\hhh({\ro 3}\aaa+{\ro 4}m)\\
 \nonumber& & +{\ro 4\left(\frac{C_7}{C_3}+1\right)}\aao\aaa\\
\label{evoluzione_g}&\leq &\Delta g-  C_3(\hhh+1)\nba+2\aao\aaa({\ro 3}m\aaa+C_8),
\end{eqnarray}
where $C_8={\ro 4m^2+2\frac{C_7}{C_3}+2}$.

\begin{Proposition}\label{GE_fin_codim}
If $\tm<\infty$, for every $\eta>0$ small enough there exists a constant $C_{\eta}>0$ depending only on $\eta$ such that the inequality
$$
\nbh\leq\eta\left| H\right|^4+C_{\eta}
$$
holds for all times.
\end{Proposition}
\proof Let $f=\nbh+ \frac{1}{C_3}(C_4+1)g-\eta\left|H\right|^4$ with $\eta>0$. By Lemma \ref{lemma_GE_codim}, Proposition \ref{evoluzione_nbh} and inequality (\ref{evoluzione_g}) we have
\begin{eqnarray*}
\dt f &\leq & \Delta f+ C_4(\hhh+1)\nba-(C_4+1)(\hhh+1)\nba\\
 & & + \frac{2}{C_3}(C_4+1)\aao\aaa({\ro 3}m\aaa+{\ro C_8}) -\eta\left(\frac{4}{m}\left|H\right|^6-12\hhh\nbh\right).
\end{eqnarray*}
We can use Lemma \ref{grad05} to find 
$$
-\left(\hhh+1\right)\nba+12\eta\hhh\nbh\leq\left(-\hhh-1+\frac{27}{4}(m+2)\eta\right)\nba,
$$
and therefore the gradient terms are non-positive for $\eta$ sufficiently small. The remaining terms are
$$
R:={\ro \frac{2}{C_3}}(C_4+1)\aao\aaa({\ro 3}m\aaa+{\ro C_8})-\frac{4\eta}{m}\left|H\right|^6.
$$
Using the pinching condition \eqref{pinching_codim_ep} we have
\begin{eqnarray*}
R&\leq& {\ro \frac{2}{C_3}}(C_4+1)\aao\left(a\hhh+b\right)\left({\ro 3}ma\hhh+C_9\right)-\frac{4\eta}{m}\left|H\right|^6,
\end{eqnarray*}
where $C_9={\ro 3}mb+{\ro C_8}$. Hence, thanks to Theorem \ref{PE_fin_codim}, we obtain
\begin{eqnarray*}
R &\leq & {\ro \frac{2}{C_3}}(C_4+1)C_0\left(\hhh+1\right)^{1-\sigma}\left(a\hhh+b\right)\left({\ro 3}ma\hhh+C_9\right)-\frac{4\eta}{m}\left|H\right|^6\\
 & \leq &  {\ro \frac{2}{C_3}}(C_4+1)C_0\left(\mu(1-\sigma)\left(\hhh+1\right)+\sigma\mu^{\frac{\sigma-1}{\sigma}}\right)\left(a\hhh+b\right)\left({\ro 3}ma\hhh+C_9\right)\\
& & -\frac{4\eta}{m}\left|H\right|^6\\ &\leq& C_{10},
\end{eqnarray*}
for some constant $C_{10}$ if $\mu$ is small enough. Putting these estimates together, we have $\dt f\leq \Delta f+C_{10}$. Since $\tm< \infty$, we conclude that there exists a constant $C_{\eta}$ depending only on $\eta$ such that $f\leq C_{\eta}$. Then, from the definition of $f$, we have
$$
\nbh\leq\nbh+{\ro \frac{1}{C_3}}(C_4+1)g\leq\eta\left|H\right|^4+ C_{\eta}.
$$
\cvd

As we have mentioned at the beginning of this section, when the codimension is greater than one we cannot repeat the proof of Proposition \ref{sez_pos}. However, using Theorem \ref{PE_fin_codim} and Proposition \ref{GE_fin_codim} we can prove that, if time is large enough, the sectional curvature of the evolving submanifold becomes positive.
\begin{Proposition}\label{k_pos_fin}
There are constants $\mu>0$ and $\vartheta>0$ such that, for any time $\vartheta<t<\tm<\infty$, the intrinsic sectional curvature of $\mm_t$ satisfies
$$
K>\mu W>0.
$$
\end{Proposition}
\proof From Gauss equation we have that
\begin{equation}\label{gauss}
2K_{ij}=2\bar K_{ij}+2\sum_{\alpha=m+1}^{2n}\left(h_{ii}^{\alpha}h_{jj}^{\alpha}-\left(h_{ij}^{\alpha}\right)^2\right),
\end{equation}
where $K_{ij}$ is the sectional curvature of $\mm_t$ of the plane spanned by two orthonormal vectors $e_i,e_j$, and $\bar K_{ij}$ is the sectional curvature of the same plane, but in $\cc\pp^n$. The idea is to use \eqref{prop_alg} restricted to the normal direction parallel to $H$. To this purpose, we fix an orthonormal basis of type (B1) with the additional requirement that $e_1,\dots,e_m$ diagonalize $h^{m+1}$, and let $\lambda_1^{m+1}\leq\dots\leq\lambda_m^{m+1}$ be the eigenvalues of $h^{m+1}$. Recalling that $\bar K\geq 1$, \eqref{gauss} becomes
\begin{eqnarray}
\nonumber 2K_{ij} & \geq & 2+2\lambda_i^{m+1}\lambda_j^{m+1}+2\sum_{\alpha=m+2}^{2n}\left(\accentset{\circ}{h}_{ii}^{\alpha}\accentset{\circ}{h}_{jj}^{\alpha}-\left(\accentset{\circ}{h}_{ij}^{\alpha}\right)^2\right)\\
\nonumber & \geq & 2+\frac{1}{m-1}\hhh-\left|h_1\right|^2-2\amo\\
\nonumber & = & 2+\frac{1}{m(m-1)}\hhh-\auo-2\amo\\
\label{gauss2} & \geq & 2+\frac{1}{m(m-1)}\hhh-2\aao.
\end{eqnarray}
By Theorem \ref{PE_fin_codim} we have
\begin{equation}\label{eq003}
2K_{ij}\geq 2+\frac{1}{m(m-1)}\hhh-2C_0\left(\hhh+1\right)^{1-\sigma}.
\end{equation}
Fix some $0<\mu<\min\{\frac{1}{2\alpha m(m-1)},\frac{1}{\beta}\}$. Then there exists $H^*$ such that, if $|H| \geq H^*$, then
\begin{equation}\label{eq004}
2+\frac{1}{m(m-1)}\hhh-2C_0\left(\hhh+1\right)^{1-\sigma}\geq 2\mu W=2\mu(\alpha\hhh+\beta).
\end{equation}
Let $\bar H(t)=\max_{\mm_t}\left|H\right|$. Since $\tm<+\infty$, we know that $\bar H(t) \to +\infty$ as $t \to \tm$, and so there
exists $\vartheta$ such that $\bar H(t) \geq \bar H$ for all $\vartheta\leq t<\tm$.
Fix some $0<\eta<\frac{1}{2}$. By Theorem \ref{GE_fin_codim}, there is a constant $C_{\eta}$ with $\left|\nb H\right|\leq\frac{1}{2}\eta^2\hhh+C_{\eta}$. By choosing larger $H^*$ and $\vartheta$ if necessary, we can assume that $C_{\eta}\leq\frac{1}{2}\eta^2(H^*) ^2$ and so $\left|\nb H\right|\leq\eta^2\bar H(t)^2$ on $\mm_t$ for $t>\vartheta$. Now fix any $t \in \,]\vartheta,\tm[$ and let $x$ be a point on $\mm_{t}$ where $\left|H\right|$ assumes its maximum. Along any geodesic starting from $x$ of length at most $r=[\eta \bar H(t)] ^{-1}$, we have $\left|H\right|\geq(1-\eta)\bar H(t)>\frac{1}{2}\bar H(t)$. By inequalities \eqref{eq003} and \eqref{eq004} we find  that
$$
K>\mu W>\mu\alpha\hhh\geq\mu\alpha\frac{\bar H(t)^2}{4}>0
$$
holds in all $B_r(x)$, with $\mu$ independent on the choice of $\eta$. Then in $B_r(x)$ we have $Ric_{ij}\geq(m-1)\frac{\mu\alpha}{4}\bar H(t)^2g_{ij}$.
Applying Myers' theorem to geodesics in $B_r(x)$ we have that, if such a geodesic has length at least
$2\pi(\bar H(t)\sqrt{\mu\alpha})^{-1}$, then it has a conjugate point. So if $\eta$ is small, precisely such that
$$
\frac{2\pi}{\bar H(t)\sqrt{\mu\alpha}}<r=\frac{1}{\eta\bar H(t)}
$$
then $B_r(x)$ covers all $\mm_t$. \cvd

To conclude the proof of the convergence of $\mm_t$ to a round point we use the main result of \cite{LXZ}, which states the following: given any Riemannian manifold with bounded geometry (in particular, the complex projective space), there is a constant $b_0>0$ such that if a submanifold of dimension $m$ satisfies 
\begin{equation}\label{lxz}
\aaa<\frac{1}{m-1}\hhh-b_0,
\end{equation}
then the mean curvature flow of this submanifold contracts to a round point in finite time. Our pinching condition \eqref{pinching_codim} on $\mm_0$ is weaker than \eqref{lxz}, but our analysis implies that  \eqref{lxz} holds on $\mm_t$ for $t$ sufficiently close to $\tm$, as the next result shows.

\begin{Proposition}\label{prop_lxz}
For every $b_0>0$, there exists a time $0<\vartheta<\tm$ such that inequality \eqref{lxz} holds on $\mm_t$  for all $\vartheta<t<\tm$.
\end{Proposition}
\proof By Theorem \ref{PE_fin_codim} we have
\begin{eqnarray*}
\aaa-\frac{1}{m-1}\hhh+b_0 & = & \aao-\frac{1}{m(m-1)}\hhh+b_0\\
 & \leq & C_0\left(\hhh+1\right)^{1-\sigma}-\frac{1}{m(m-1)}\hhh+b_0,
\end{eqnarray*}
which is negative at the points $(x,t)$ where $\hhh(x,t)$ is big enough. Using Myers' theorem as in the the proof of Proposition \ref{k_pos_fin} we obtain the assertion. \cvd

\section{Infinite maximal time}\setcounter{equation}{0}\setcounter{Theorem}{0}
Throughout this section we assume $\tm=\infty$. In this case, the argument is simpler than in the case of finite maximal time, because the improvement of pinching can be obtained directly from the maximum principle, as shown in the next result.

\begin{Proposition}\label{PE_inf_codim}
There are positive constants $C_0$ and $\delta_0$ depending only on the initial manifold $\mm_0$ such that
\begin{equation*}
\aao\leq C_0\left(\hhh+1\right)e^{-\delta_0t}
\end{equation*}
holds for any time $0\leq t<\tm=\infty$.
\end{Proposition}
\proof Using Proposition \ref{fs_codim} with $\sigma=0$ and the maximum principle, we have that
$$
f_0\leq C_0'e^{-\delta_0 t},
$$
for some positive constants $C_0'$ and $\delta_0$ that depend only on the initial data. Recalling that 
$$
f_0=\frac{\aao}{\alpha\hhh+\beta},
$$
we obtain the assertion for an appropriate constant $C_0$.\cvd

Note that the above result is trivial for small values of $t$, while it becomes significant when $t$ is arbitrarily large. As a first consequence of this estimate, we can prove that the intrinsic sectional curvature of the evolving submanifold becomes positive for time large enough, similarly to the case of finite maximal time.

\begin{Proposition}\label{k_pos_inf}
There are constants $\mu>0$ and $\vartheta>0$ such that, for any time $\vartheta<t<\tm=\infty$, the intrinsic sectional curvature of $\mm_t$ satisfies
$$
K>\mu W>0.
$$
\end{Proposition}
\proof As in the proof of Proposition \ref{k_pos_fin}, we have $2K_{ij}\geq 2+\frac{1}{m(m-1)}\hhh-2\aao.$
By the exponential decay of $\aao$ proved in Proposition \ref{PE_inf_codim}, we have
$$
2K_{ij}  \geq  2+\frac{1}{m(m-1)}\hhh-2C_0\left(\hhh+1\right)e^{-\delta_0 t}
  \geq  2\mu W>0,
$$
for $\mu>0$ small enough and $t$ sufficiently big. \cvd

Now we can follow a procedure similar to the previous section.

\begin{Lemma}\label{H4inf} 
There exists $C_7>0$ such that
\begin{eqnarray*}\dt\hhh\aao &\leq&  \Delta(\hhh\aao)- {C_3}\hhh\nba+C_7\nba\\
&&+2\hhh\aao({\ro 3}\aaa+{\ro 4}m).
\end{eqnarray*}
\end{Lemma} 
\proof We proceed like in the proof of Lemma \ref{H4fin}, but this time we use Proposition \ref{PE_inf_codim}, to find that
\begin{eqnarray*}
{\ro 6}\left|H\right|\sqrt{m+2}\nba\aao &\leq &{\ro 6}\left|H\right|\sqrt{m+2}\nba\sqrt{C_0(\hhh+1)}e^{-\delta_0t/2}\\
 & \leq &  {C_3}\hhh\nba+C_7\nba
\end{eqnarray*}
if $C_7$ is chosen large enough. In fact, for $t$ large enough this follows from the exponential decay of $\aao$, while if $t$ varies on any compact interval of $[0,\infty[$ it follows from the boundedness of $\hhh$. \cvd

Now we consider the function $g$ defined in \eqref{defg}. Using Lemma \ref{lemma_GE_codim} and Lemma \ref{H4inf}, we can repeat the computations of the previous sections to conclude that inequality \eqref{evoluzione_g} holds also in this case. We can now prove a gradient estimate for the curvature.

\begin{Theorem}\label{GE_inf_codim}
For every $\eta>0$ small enough there exists a constant $C_{\eta}>0$ depending only on $\eta$ such that for all time we have the estimate
$$
\nbh\leq\left(\eta\left| H\right|^4+C_{\eta}\right)e^{-\delta_0 t/2}.
$$
\end{Theorem}
\proof The proof is similar to Proposition \ref{GE_fin_codim}. Let us define $$f=e^{\delta_0 t/2}\left(\nbh+ \frac{1}{C_3}(C_4+\delta_0m)g\right)-\eta\left|H\right|^4.$$  By Proposition \ref{evoluzione_nbh}, Lemma \ref{lemma_GE_codim} and inequality (\ref{evoluzione_g}) we have
\begin{eqnarray*}
\dt f &\leq &\Delta f+\left[\frac{\delta_0}{2}\nbh+ \frac{\delta_0}{2C_3}\left(C_4+\delta_0m\right)\left(\hhh\aao+2(C_7+1)\aao\right)\right]e^{\delta_0 t/2}\\
 & & +\left[-\delta_0 m(\hhh+1)\nba+ \frac{2}{C_3}(C_4+\delta_0m)\aao\aaa({\ro 3}m\aaa+C_8)\right]e^{\delta_0 t/2}\\
 & & -\eta\left(\frac{2}{m}\left|H\right|^6-12\hhh\nbh\right).
\end{eqnarray*}
By Lemma \ref{grad05}, the gradient terms satisfy
$$
\begin{array}{l}
\displaystyle{\left[\frac{\delta_0}{2}\nbh-\delta_0 m\left(\hhh+1\right)\nba\right]e^{\delta_0 t/2}+12\eta\hhh\nbh }\\
\displaystyle{\qquad\leq\left[\frac{\delta_0}{2}-\frac{16\delta_0 m}{9(m+2)}(\hhh+1)+12\eta\hhh\right]\nbh e^{\delta_0 t/2}  }
\end{array}
$$
and therefore they are non-positive for $\eta$ sufficiently small. 
We call $R$ the remaining terms. Using condition \eqref{pinching_codim} and Theorem \ref{PE_inf_codim}, we can find a constant $\Lambda$ such that
\begin{eqnarray*}
R&\leq& C_0\Lambda\left(\hhh+1\right)\left(\left|H\right|^4+1\right)e^{-\delta_0 t/2}-\frac{2\eta}{m}\left|H\right|^6\\
 &\leq&\left[ C_0\Lambda\left(\hhh+1\right)\left(\left|H\right|^4+1\right)e^{-\delta_0 t/4}-\frac{2\eta}{m}\left|H\right|^6\right]e^{-\delta_0 t/4}\\
 &\leq &C_{10} e^{-\delta_0 t/4},
\end{eqnarray*}
for a suitably large constant $C_{10}$. Note that this is true, because $e^{-\delta_0 t/4}$ is small, for $t$ big enough, and because $\hhh$ is bounded, for $t$ small. Then there exists a constant $C_{\eta}$ such that $f\leq C_{\eta}$. Recalling the definition of $f$ we conclude the proof.\cvd

We now show that, if $\tm=\infty$, the curvature is uniformly bounded.

\begin{Lemma}\label{nosing}
If $\tm=\infty$, then $\hhh$ is bounded uniformly for all $t$.
\end{Lemma}
\proof Let $b_0$ the constant which appears in the main theorem in \cite{LXZ}. From Proposition \ref{PE_inf_codim} we have
\begin{eqnarray*}
\aaa-\frac{1}{m-1}\hhh+b_0 & = & \aao-\frac{1}{m(m-1)}\hhh+b_0\\
 & 	\leq & C_0\left(\hhh+1\right)e^{-\delta_0 t}-\frac{1}{m(m-1)}\hhh+b_0.
\end{eqnarray*}
Observe that the right-hand side is negative if $t$ and $\hhh$ are big enough. Using Proposition \ref{k_pos_inf} and Theorem \ref{GE_inf_codim}, we can apply Myers' theorem like in the proof of Proposition \ref{prop_lxz} to show that, if the curvature is sufficiently large at some point, then it is large everywhere. Therefore, if $\hhh$ becomes arbitrarily large as $t \to \infty$ we obtain that, for $t$ big enough, 
$$
\aaa-\frac{1}{m-1}\hhh+b_0<0
$$
everywhere on $\mm_t$. Then the main theorem in \cite{LXZ}  implies that the mean curvature flow with initial value $\mm_t$ shrinks to a point in finite time, giving a contradiction. \cvd

Now we have all the ingredients to prove the convergence in the case $\tm=\infty$. Since $\hhh$ stays bounded, Proposition \ref{PE_inf_codim} and Theorem \ref{GE_inf_codim} give that there is a constant $C$ such that
$$
\aao\leq Ce^{-\delta_0 t}, \qquad \nbh\leq Ce^{-\delta_0 t/2}.
$$
Applying once again Myers' theorem, the diameter of $\mm_t$ is uniformly bounded and so $\hhh_{max}-\hhh_{min}\leq Ce^{-\delta_0 t/2}$. Moreover $\hhh_{min}=0$ otherwise the evolution equation for $\hhh$ from Lemma \ref{evoluzione_AH}, together with \eqref{e.R2}, would imply by a standard comparison argument the finite time blow up of $\hhh$, in contradiction with the assumption $\tm=+\infty$. Then $\hhh$ decays exponentially fast and 
$$
\aaa=\aao+\frac{1}{m}\hhh\leq Ce^{-\delta_0 t/2},
$$
for some $C>0$. We deduce
$$
\int_0^{\infty}\left|\dt g_{ij}\right|dt=\int_0^{\infty}\left|H\right|\left|A\right|dt\leq \sqrt{m}\int_0^{\infty}\aaa dt\leq \sqrt{m}C\int_0^{\infty}e^{-\delta_0t/2}\leq \bar C,
$$
for some $\bar C>0$. So we can apply a result by Hamilton \cite[Lemma 14.2]{Ha} to obtain that there is a continuous limit metric $g_{ij}(\infty)$. By the same method used in \cite[\S 10]{H1}, we can show that the exponential decay for $\aaa$ gives the exponential decay for all derivatives $\nb^k A$ by means of interpolation inequalities. This finally gives $C^{\infty}$-convergence to a smooth totally geodesic submanifold $\mm_{\infty}$. By our smallness assumption on the codimension $k$, the only possibility is that $\mm_{\infty}=\cc\pp^{n'}$ for some $n'<n$ as implied by Theorem 3.25 of \cite{Be1}. Therefore, if $k$ is odd this possibility cannot happen and we can only have a singularity in finite time.
This concludes the proof of Theorem \ref{maincodim}.  

\section{Extensions to quaternionic projective spaces}\setcounter{equation}{0}\setcounter{Theorem}{0}
In this last section we show that in the case of hypersurfaces our main result, Theorem \ref{maincodim}, can be easily extended to the flow in a quaternionic projective space. Let $\kk$ be either the field $\cc$ of complex numbers or the associative algebra $\qq$ of quaternions, and let $c$ be a positive constant. We denote by $\kk\pp^n(4c)$ the projective space over $\kk$ with sectional curvature $c\leq\bar K\leq 4c$, and we consider the mean curvature flow of a real hypersurface of  $\kk\pp^n(4c)$.

\begin{Theorem}\label{mainCROSSes}
Let $n\geq 3$, $c>0$, and let $\mm_0$ be a closed real hypersurface of $\kk\pp^n(4c)$. Let $m$ be the real dimension of $\mm_0$ and suppose that $\mm_0$ satisfies
\begin{equation}\label{pinching_CROSSes}
\aaa<\frac{1}{m-1}\hhh+2c.
\end{equation}
Then the mean curvature flow with initial condition $\mm_0$ has a smooth solution $\mm_t$ on a finite time interval $0\leq t<\tm<\infty$ and the flow converges to a round point as $t$ goes to $\tm$.
\end{Theorem}
The proof is the same exposed in the previous sections for the case of hypersurfaces of $\cc\pp^n=\cc\pp^n(4)$. The constants used are
\begin{equation*}
m=\left\{\begin{array}{ll}
2n-1 & \text{if }\kk=\cc,\medskip \\
4n-1 & \text{if }\kk=\qq,
\end{array}
\right.
\qquad\text{and}\qquad
\rr=\left\{\begin{array}{rl}
2(n+1)c & \text{if }\kk=\cc,\medskip \\
4(n+2)c & \text{if }\kk=\qq.
\end{array}\right.
\end{equation*}
As we have observed in the complex case, the proof that the flow develops a singularity in finite time is in some sense indirect and is related to the global structure of the projective spaces we are considering. Namely, we show that a solution defined for all times would converge to a totally geodesic hypersurface, but this is excluded because in $\kk\pp^n(4c)$ there are no such hypersurfaces.

Theorem \ref{mainCROSSes} implies the following classification result.
\begin{Corollary}
Let $n\geq 3$ and $c>0$.
\begin{enumerate}
\item If $\mm_0$ is a closed real hypersurface of $\kk\pp^n(4c)$ satisfying the pinching condition \eqref{pinching_CROSSes}, then $\mm_0$ is diffeomorphic to a sphere.
\item For any minimal closed real hypersurface of $\kk\pp^n(4c)$, $\aaa\geq 2c$ holds.
\end{enumerate}
\end{Corollary}

\noindent Theorem \ref{mainCROSSes} is the generalization of the main theorem of \cite{H3} about pinched hypersurfaces of the sphere to all CROSSes (compact rank-one symmetric spaces) with sufficiently large dimension. Unfortunately, these techniques do not allow to find an analogous result for the Cayley plane $\cc a\pp^2$.
 
The next example shows that Theorem \ref{mainCROSSes} is not a trivial consequence of the general result in \cite{H2}, because there are non-convex hypersurfaces in the class considered.

\begin{Example}{\rm
Consider for simplicity $c=1$ and let $\mm_0$ be a geodesic sphere in $\cp^n$. In \cite{NR} it is proved that $\mm_0$ has two distinct principal curvatures: $\lambda_1=2\cot(2u)$ with multiplicity $1$ and $\lambda_2=\cot(u)$ with multiplicity $2(n-1)$, for some $0< u <\frac{\pi}{2}$. For any $u>\frac{\pi}{4}$, we have $\lambda_1<0$ and $\lambda_2>0$, so $\mm_0$ is not convex. Moreover, it is easy to compute that in this case condition \eqref{pinching_CROSSes} is equivalent to
$$
2(2n-3)\cot^2(2u)-2(n-1)\cot^2(u)<0.
$$
Hence, there are non-convex examples in our class for every $n$.  In the same way, a geodesic sphere in $\qq\pp^n$ has principal curvatures $\lambda_1=2\cot(2u)$ with multiplicity $3$ and $\lambda_2=\cot(u)$ with multiplicity $4(n-1)$, for some $0< u <\frac{\pi}{2}$ (see for example \cite{MP}). Condition \eqref{pinching_CROSSes} in this case becomes
$$
3(4n-5)\cot^2(2u)-4(n-1)\cot^2(u)+4n-5<0,
$$
so we have non-convex examples in our class for $\kk=\qq$ too. We remark that, even if the initial hypersurface is not convex, it becomes convex approaching the maximal time, as a consequence of the convergence to a round point.
}\end{Example}
\medskip

\noindent {\bf Acknowledgments} 
The results of this paper are part of Giuseppe Pipoli's PhD thesis, written at the Department of Mathematics, University ``Sapienza'' of Rome. 
Giuseppe Pipoli was partially supported by PRIN07 ``Geometria Riemanniana e strutture differenziabili'' of MIUR (Italy) and Progetto universitario Univ. La Sapienza ``Geometria differenziale -- Applicazioni''. Carlo Sinestrari was partially supported by FIRB--IDEAS project ``Analysis and beyond'' and by the group GNAMPA of INdAM (Istituto Nazionale di Alta Matematica).

\bigskip

\noindent Giuseppe Pipoli, Institut Fourier, Universit\'e Joseph Fourier (Grenoble I), UMR 5582, CNRS-UJF, 38402, Saint-Martin-d'H\`eres, France. E-mail: giuseppe.pipoli@ujf-grenoble.fr \\

\noindent Carlo Sinestrari, Dipartimento di Matematica, Universit\`a di Roma ``Tor Vergata'', Via della Ricerca Scientifica, 00133, Roma, Italy. E-mail: sinestra@mat.uniroma2.it

\end{document}